\RequirePackage{fix-cm}
\documentclass[smallcondensed]{svjour3}     
\smartqed  
\usepackage{graphicx}
\usepackage{amsmath,amssymb}
\usepackage{bbm}
\usepackage{hyperref}
\usepackage{graphicx}
\usepackage{subfig}
\usepackage{epstopdf}
\usepackage{bbold}	  
\usepackage{algpseudocode,algorithm,algorithmicx}
\usepackage{cite}
\begin{document}

\title{A fast sparse spectral method for nonlinear integro-differential Volterra equations with general kernels
}

\titlerunning{A sparse spectral method for nonlinear integro-differential Volterra equations}        

\author{Timon S. Gutleb
}


\institute{Timon S. Gutleb \at
              Department of Mathematics, Imperial College London, UK \\
              \email{t.gutleb18@imperial.ac.uk}           
}

\maketitle

\begin{abstract}
We present a sparse spectral method for nonlinear integro-differential Volterra equations based on the Volterra operator's banded sparsity structure when acting on specific Jacobi polynomial bases. The method is not restricted to convolution-type kernels of the form $K(x,y)=K(x-y)$ but instead works for general kernels at competitive speeds and with exponential convergence. We provide various numerical experiments on problems with or without known analytic solutions and comparisons with other methods.
\end{abstract}

\section{Introduction}
The \emph{Volterra integral operator} is defined by
\begin{equation}\label{eq:volterraintegral}
    (\mathcal{V}_K u)(x) := \int_0^{x} K(x,y) u(y) \mathrm{d}y.
\end{equation}
Integral equations involving Volterra operators occur in diverse applications in various disciplines of engineering \cite{unterreiter_volterra_1996,zakes_application_2016}, finance and economics \cite{apartsyn_classes_2014,nedaiasl_product_2019} as well as the natural sciences \cite{micke_application_1990,hethcote_integral_1980,geiser_iterative_2013,pruss_evolutionary_2012, van_den_bosch_pandemics_1999,krimer_non-markovian_2014,krimer_sustained_2016,xiang_efficient_2013}. Beyond \emph{linear} Volterra integral equations of first and second kind, which respectively take the forms
\begin{equation*}
\mathcal{V}_K u = g \qquad\hbox{or}\qquad  ( \lambda I + \mathcal{V}_K) u = g,
\end{equation*}
there exist a vast range of possible generalizations, the most important of which are linear Volterra integro-differential equations (VIDEs), nonlinear Volterra integral equations and in particular nonlinear VIDEs, the last of which take the form:
\begin{align}\label{eq:generalcase}
\sum_{k=0}^m \lambda_k \frac{\mathrm{d}^k}{\mathrm{d}x^k}u(x) = g + \mathcal{V}_K f(u),
\end{align}
where $m \in \mathbb{N}$ and $\forall i: \lambda_i \in \mathbb{R}$. Linear Volterra equations thus correspond to the special case $f(u)=u$. For certain very well-behaved cases one can use variational iteration, Adomian decomposition or Laplace transform methods \cite{wazwaz_linear_2011} to try to obtain analytic solutions for such equations but in general one must use numerical methods to find approximate solutions. We assume in this paper that the equations we intend to solve with our proposed numerical scheme are solvable with unique solutions and omit discussion of ill-posed problems. A number of results on criteria for the existence of solutions are known, see e.g. \cite{meehan_existence_1998,gordji_existence_2011,zhang_existence_2018} and the references therein. Interest in efficient algorithms with good convergence properties is high, resulting in a variety of competitive methods for linear, nonlinear and integro-differential Volterra equations. For decades \cite{brunner_numerical_1973} researchers have been proposing various forms of increasingly refined collocation and iteration methods to approach nonlinear VIDEs \cite{song_analysis_2019, allaei_collocation_2017,driscoll_automatic_2010,agbolade_solutions_2017,brunner_collocation_2004,shayanfard_collocation_2019}. More recently they were also joined by a number of wavelet-based methods \cite{saeedi_numerical_2011,biazar_chebyshev_2012,zhu_solving_2012,sahu_legendre_2015,lepik_haar_2006,heydari_wavelets_2014}. Numerical solvers based on homotopy perturbation methods have also been proposed \cite{ghasemi_numerical_2007}. A highly efficient spectral solver for the case of convolution kernels $K(x,y)=K(x-y)$ for linear VIDEs using low rank approximations was recently described by Hale \cite{hale_ultraspherical_2019}. Gutleb and Olver described a general kernel sparse spectral method for linear Volterra integral equations in \cite{gutleb_sparse_2019}, which forms the background of the present paper.\\
In this paper we present a spectral method which works for linear, nonlinear, integro-differential and many other types of Volterra equations of first, second and third kind while utilizing polynomial spaces in which the Volterra operators have a sparse banded structure. As such this paper is a direct generalization of results in \cite{gutleb_sparse_2019}, which derived said banded structure of the Volterra operator in Jacobi polynomial bases and on the basis of \cite{olver_sparse_2019} developed an efficient Clenshaw algorithm based approach to numerically generate the operator. This method thus combines the unmatched exponential convergence of spectral methods with highly efficient sparse linear algebra, a very promising combination which has been successful in recent years \cite{snowball_sparse_2019,hale_fast_2018,townsend_automatic_2015,olver_fast_2013,olver_practical_2014}. In addition to efficiency via bandedness and exponential convergence rate, the proposed method is furthermore not limited to convolution kernel cases, i.e. kernels of form $K(x,y)=K(x-y)$, a common restriction in competitively fast and accurate approaches \cite{hale_ultraspherical_2019}. Due to the extensive literature and code libraries available on numerical solutions for linear, nonlinear VIDEs an exhaustive comparison with other methods is far beyond the scope of a single paper. We will however present comparisons with recent collocation methods, as these are the most competitive and ubiquitous solvers available.\\

The structure of this paper is as follows: In sections \ref{sec:functionapproximationuni} and \ref{sec:functionapproximationmulti} we introduce the necessary framework of uni- and multivariate orthogonal polynomial approximation of functions. Section \ref{sec:volterrabanded} briefly recounts relevant results from \cite{gutleb_sparse_2019}, detailing the banded structure of the Volterra operator on appropriately chosen triangle domain Jacobi polynomial bases. Section \ref{sec:integrodifferential} details the extension of the linear Volterra integral equation method to a general linear integro-differential case by augmenting the system with appropriate evaluation operators. Section \ref{sec:nonlinear} details the extension to nonlinear Volterra integral equations using an iterative approach and describes further relevant numerical experiments. Section \ref{sec:nonlinearintegrodiff} combines these two approaches, finally extending the method to general kernel nonlinear VIDEs. Section \ref{sec:numexprsection} showcases various numerical experiments for linear VIDEs as well as nonlinear VIEs and VIDEs to verify and test applicability, convergence rate and competitiveness including comparisons to collocation methods in Chebfun. We close with notes on convergence for the methods proposed in this paper in section \ref{sec:analysis} and discuss applicability and potential further research directions in the conclusion.
\subsection{Function approximation with univariate orthogonal polynomials} \label{sec:functionapproximationuni}
In what follows we introduce the relevant elements of function approximation in univariate orthogonal polynomial bases, primarily focusing on the Jacobi polynomials, which are required to understand the proposed spectral method for nonlinear VIDEs. The reasons for highlighting the specific chosen properties above others will become clear when we discuss the methods in sections \ref{sec:integrodifferential} and \ref{sec:nonlinear}. For a more general and complete introduction into the theory of function approximation with univariate orthogonal polynomials, see e.g. \cite{gautschi_orthogonal_2004,beals_special_2016}.\\
The Jacobi polynomials $P_n^{(\alpha, \beta)}(t)$ are a univariate complete set of polynomials orthogonal with respect to the weight function $(1-t)^\alpha (1+t)^\beta$ on their natural domain $t \in [-1,1]$, meaning they satisfy
\begin{equation*}
\int_{-1}^1 (1-x)^{\alpha} (1+t)^{\beta} P_n^{(\alpha,\beta)} (t)P_m^{(\alpha,\beta)} (t)\,dt =\tfrac{2^{\alpha+\beta+1}}{2n+\alpha+\beta+1} \tfrac{\Gamma(n+\alpha+1)\Gamma(n+\beta+1)}{n!\Gamma(n+\alpha+\beta+1)} \delta_{nm},
\end{equation*}
where $\alpha, \beta \geq -1$. As such, the Legendre polynomials correspond to the special case $\alpha = \beta = 0$ and the ultraspherical or Gegenbauer polynomials correspond to the special case in which $\alpha = \beta$.
Spectral methods can make use of the fact that complete sets of orthogonal polynomials can be used to approximate any sufficiently smooth function $f(t)$ defined on a real interval $\left(a,b \right)$ via the expansion
\begin{equation*}
    f(t) = \sum_{n=0}^\infty p_n(x) f_n = \mathbf{p}(t)^\mathsf{T} \mathbf{f},
\end{equation*}
where $f_n$ are the unique constant coefficients of $f(t)$ in the given complete polynomial basis $\mathbf{p}(t)$ which is orthogonal on the domain $(a,b)$. These coefficients $f_n$ may be computed efficiently for various sets of orthogonal polynomials using methods and C libraries by Slevinsky \cite{slevinsky_conquering_2017,slevinsky_fast_2017,slevinsky_fasttransforms_2019}, while evaluation of polynomials is efficiently performed using Clenshaw's algorithm, see e.g. \cite{gautschi_orthogonal_2004}. While the interval $\left[-1,1 \right]$ is the natural choice for the Jacobi polynomials one can easily shift them to any other real interval required by an application. The method in this paper exclusively makes use of the Jacobi polynomials shifted to the unit interval $[0,1]$, so we introduce the following shorthand notation:
\begin{align*}
\mathbf{\tilde{P}}(x) = \mathbf{P}(2x-1), \quad x\in[0,1].
\end{align*}
Once expanded in the above way, performing addition and subtraction of functions has an obvious element-wise implementation. Additionally, being orthogonal polynomials the Jacobi polynomials satisfy a three term recurrence relationship
\begin{equation*}
t P_{n}^{(\alpha,\beta)}(t) = c_{n-1} P_{n-1}^{(\alpha,\beta)}(t) + a_n  P_{n}^{(\alpha,\beta)}(t) + b_n P_{n+1}^{(\alpha,\beta)}(t), \quad n\geq1,
\end{equation*}
 which allows for efficient computation of function multiplication in this framework via a tridiagonal so-called Jacobi operator:
\begin{align*}
    &t f(t) = \mathbf{P}(t)^\mathsf{T} \mathrm{J}^\mathsf{T} \mathbf{f},\\
&\mathrm{J} = \begin{pmatrix} a_0 & b_0 & &  \\ c_0 & a_1 & b_1 & \\ & c_1 & a_2 & \ddots \\ &&\ddots&\ddots \end{pmatrix}.
\end{align*}
The exact elements of the Jacobi operator depend on the Jacobi parameters $(\alpha, \beta)$ of the chosen basis, see e.g. \cite[18.9(i)]{nist_2018} for explicit values of $a_i, b_i$ and $c_i$. As the sparsity of the operators in this paper relies heavily on correctly moving between bases with different Jacobi parameters $(\alpha, \beta)$ we will index coefficient vectors with the Jacobi parameters of the basis of expansion, i.e. by writing $\mathbf{f}_{(\alpha,\beta)}$, where it might otherwise be ambiguous. The above properties allow for the development of software packages capable of performing arithmetic on functions using highly efficient sparse linear algebra, where functions are replaced by coefficient vectors. One such package is ApproxFun.jl \cite{olver_juliaapproximation2019} writen in the Julia programming language \cite{beks2017}, which is used as the background environment of the implementations presented in this paper. Other available software packages include among others the Dedalus project \cite{burns_dedalus:_2019}, Chebfun \cite{pachon_piecewise-smooth_2010,battles_extension_2004,driscoll_chebop_2008}.\\
Beyond arithmetic with functions, there are a number of other useful operators in the sparse spectral method toolbox. We may, for example, freely shift from a basis with Jacobi parameters $(\alpha,\beta)$ to one with higher parameters $(\alpha+n,\beta+m)$ via
\begin{align*}
    \mathbf{\tilde{P}}^{(\alpha+n, \beta+m)}(x)^\mathsf{T} \mathrm{S}_{(\alpha, \beta)}^{(\alpha+n, \beta+m)} \mathbf{f}_{(\alpha,\beta)} = \mathbf{\tilde{P}}^{(\alpha+n, \beta+m)}(x)^\mathsf{T} \mathbf{f}_{(\alpha+n, \beta+m)},
\end{align*}
where $n,m \in \mathbb{N}$ using a sequence of upper bidiagonal raising operators \cite[18.9.5]{nist_2018}
\begin{equation*}
\mathrm{S}_{(\alpha, \beta)}^{(\alpha+n, \beta)} = \mathrm{S}_{(\alpha+n-1, \beta)}^{(\alpha+n, \beta)}   \cdots \mathrm{S}_{(\alpha+1, \beta)}^{(\alpha+2, \beta)} \mathrm{S}_{(\alpha, \beta)}^{(\alpha+1, \beta)}
\end{equation*}
and an analogous sequence of operators for the second Jacobi parameter. Lowering operators are also available, although this is in general only possible in a sparse (lower bidiagonal) way with added weights \cite[18.9.6]{nist_2018}:
\begin{align*}
    x f(x) = \mathbf{\tilde{P}}^{(\alpha-1, \beta)}(x)^\mathsf{T} \mathrm{L}_{(\alpha,\beta)}^{(\alpha-1,\beta)} \mathbf{f}_{(\alpha,\beta)},\\
    (1-x) f(x) = \mathbf{\tilde{P}}^{(\alpha, \beta-1)}(x)^\mathsf{T} \mathrm{L}_{(\alpha,\beta)}^{(\alpha,\beta-1)} \mathbf{f}_{(\alpha,\beta)}.
\end{align*}
We may also mirror functions on their domain in a given Jacobi polynomial basis by using the very useful symmetry property \cite[Table 18.6.1]{nist_2018}:
\begin{align} \label{eq:symmetryJacobi}
P_n^{(\alpha,\beta)}(-x)= (-1)^n P^{(\beta,\alpha)}_n(x).
\end{align}
In particular, on  $[0,1]$ we can define a diagonal reflection operator via
\begin{align*}
   f(1-x) = \sum_n (-1)^n f_{(\alpha,\beta),n} \tilde{P}_n^{(\beta, \alpha)}  = \mathbf{\tilde{P}}^{(\alpha, \beta)}(x)^\mathsf{T} \mathrm{R} \mathbf{f}_{(\alpha,\beta)} .
\end{align*}
Differentiation is also a sparse (diagonal) operation if we simultanously increment the Jacobi parameters \cite[18.9.15]{nist_2018}, i.e.:
\begin{align}\label{eq:derivativeop1}
\frac{\mathrm{d}}{\mathrm{d}x} f(x) &=  \sum_n f_{(\alpha,\beta),n} \frac{\mathrm{d}}{\mathrm{d}x} P_n^{(\alpha,\beta)}(x)\\ &= \sum_n f_{(\alpha,\beta),n} \frac{1}{2} \left( n+\alpha+\beta+1 \right) P_{n-1}^{(\alpha+1,\beta+1)}(x)\\ &= \mathbf{\tilde{P}}^{(\alpha+1,\beta+1)}{}(x)^\mathsf{T} \mathcal{D}^{(\alpha, \beta)} \mathbf{f}_{(\alpha,\beta)}.\label{eq:derivativeop3}
\end{align}
Importantly, this means that repeated sparse differentiation is not equivalent to a repeat application of the same operator $\mathcal{D}^{(\alpha,\beta)}$. As the derivative operator shifts coefficient vectors to a higher parameter basis the second derivative operator is actually a combination of two distinct derivative operators acting on different bases and so on for higher derivatives. We thus denote the $n$-th derivative operator acting on a coefficient vector in $\mathbf{\tilde{P}}^{(\alpha,\beta)}$ basis by $ \mathcal{D}^{(\alpha,\beta)}_n$, where
\begin{align*}
\frac{\mathrm{d}^n}{\mathrm{d}x^n} f(x) &= \mathbf{\tilde{P}}^{(\alpha+n,\beta+n)}{}(x)^\mathsf{T} \mathcal{D}^{(\alpha+n-1,\beta+n-1)} \cdots  \mathcal{D}^{(\alpha,\beta)} \mathbf{f}_{(\alpha,\beta)}\\ &= \mathbf{\tilde{P}}^{(\alpha+n,\beta+n)}{}(x)^\mathsf{T}  \mathcal{D}^{(\alpha,\beta)}_n \mathbf{f}_{(\alpha,\beta)},
\end{align*}
instead of the potentially misleading notation $\mathcal{D}^n$ which may evoke false intuitions of commutativity of the operators. The last component of theory we need for our univariate function approximation purposes are endpoint evaluation operators which will be used to enforce boundary conditions in integro-differential equations. From the viewpoint described above, functions are coefficient vectors and multiplications, derivatives and basis changes are operators on coefficient vectors (matrices in finite dimensional approximation space). Functionals, e.g. evaluation operators $\mathcal{E}$ at an endpoint, must act on coefficient vectors to return a scalar value and are thus represented by row vectors. In particular, for the Jacobi polynomials we can make use of the known property \cite[Table 18.6.1]{nist_2018}:
\begin{align*}
f(1) = \mathbf{P}^{(\alpha,\beta)}{}(x)^\mathsf{T} \mathcal{E}_1 \mathbf{f}_{(\alpha,\beta)} = \sum_n f_{(\alpha,\beta),n} P_n^{(\alpha,\beta)}(1) = \sum_n f_{(\alpha,\beta),n}  \frac{(\alpha + 1)_n}{n!},
\end{align*}
where $(\cdot)_n$ denotes the Pochhammer symbol or rising factorial \cite[5.2(iii)]{nist_2018}. Via the symmetry property in Equation (\ref{eq:symmetryJacobi}) we obtain a similar evaluation operator for the other endpoint of our chosen interval domain.

\subsection{Function approximation with multivariate orth. polynomials} \label{sec:functionapproximationmulti}
This section introduces required elements of function approximation in multivariate orthogonal polynomial bases, focusing on the Jacobi polynomials on the triangle domain $$T^2 = \{(x,y):0\leq x \leq 1, 0 \leq y \leq 1-x \},$$
which was discussed in detail in \cite{olver_sparse_2019}. Multivariate polynomials are not yet widely used in numerical methods despite their vast potential. For a more general and complete introduction to theoretical aspects of multivariate orthogonal polynomials, see \cite{dunkl_orthogonal_2014}.\\
Function approximation with multivariate orthogonal polynomials works in direct analogy to the univariate case. For a given multivariate function $f(x,y)$ defined on some suitable $2$-dimensional domain and given a set of complete multivariate orthogonal polynomials on said domain, we may expand the function via
\begin{align*}
f(x,y) = \sum_{n=0}^\infty \sum_{k=0}^n f_{nk} p_{nk}(x,y) = \mathbf{p}(x,y){}^\mathsf{T} \mathbf{f}.
\end{align*}
A generalization to $n$-dimensional cases is straightforward, cf.  \cite{dunkl_orthogonal_2014}. On the triangle $T^2$, a sensible choice of polynomial basis is found in the triangle Jacobi polynomials, also known as Proriol polynomials, which are defined via reference to the univariate Jacobi polynomials \cite[Proposition 2.4.1]{dunkl_orthogonal_2014}:
\begin{equation*}
    P_{k,n}^{(\alpha, \beta, \gamma)} (x,y) = (1-x)^k \tilde{P}_{n-k}^{\left(2k+\beta+\gamma+1,\alpha\right)}\left(x\right) \tilde{P}_k^{\left(\gamma,\beta\right)}\left(\frac{y}{1-x}\right).
\end{equation*}
Expansion coefficients for functions on the triangle Jacobi polynomial basis, such as required for the kernel $K(x,y)$ discussed in the next section, may be computed efficiently using C libraries by Slevinsky \cite{slevinsky_conquering_2017,slevinsky_fast_2017,slevinsky_fasttransforms_2019}. As in the univariate case, we can define a variety of operators acting on coefficient space such as multiplication operators based on Jacobi operators, derivative operators and basis change operators \cite{olver_sparse_2019}. The novelty is that for 2-dimensional spaces such as the triangle we need to distinguish between the $x$ and $y$ variables and thus have two different Jacobi operators: $\mathrm{J}_x$ for the $x$ variable and $\mathrm{J}_y$ for the $y$ variable, which are now block tri-diagonal operators instead of being tridiagonal, see \cite{olver_sparse_2019} for details.
\subsection{Banded sparsity of the linear Volterra operator in Jacobi bases} \label{sec:volterrabanded}
It was shown in \cite{gutleb_sparse_2019} that the Volterra integral operator is sparse with banded structure on appropriate Jacobi polynomial spaces. Based on this, a sparse spectral method with exponential convergence for linear Volterra equations with general kernels was motivated and analyzed. The results are based on interpreting the Volterra operator as acting on multivariate Jacobi bases on a triangle domain. The idea behind the linear method follows the schemes in Algorithms \ref{alg:linear1} and \ref{alg:linear2}. In this section we briefly review these methods to the degree necessary to follow the integro-differential and nonlinear extension in this paper. For the full discussion of the linear case, we refer to \cite{gutleb_sparse_2019}.\\
The move to the triangle domain may initially be motivated by noting that the Volterra integral operator $\int_0^{l(x)} K(x,y) u(y) \mathrm{d}y$ acting on $u$ may be considered for other upper bounds $l(x)$ than $x$, in particular $l(x)=1-x$. The Proriol polynomials with parameters $(0,0,0)$, being orthogonal on the triangle domain, behave well with respect to this integration:
\begin{align*}
\int_{0}^{1-x} f(x,y) \mathrm{d}y &= \int_{0}^{1-x} \sum_{n=0}^\infty \sum_{k=0}^n p_{n,k}(x,y) f_{n,k} \mathrm{d}y  \\ &= \sum_{n=0}^\infty \sum_{k=0}^n f_{n,k} (1-x)^k \tilde{P}_{n-k}^{(2k+1,0)}(x) \int_{0}^{1-x} \tilde{P}_{k}^{(0,0)}\left(\frac{y}{1-x}\right) \mathrm{d}y \\ &= \sum_{n=0}^\infty \sum_{k=0}^n f_{n,k} (1-x)^{k+1} \tilde{P}_{n-k}^{(2k+1,0)}(x) \int_{0}^{1} \tilde{P}_{k}^{(0,0)}\left(s\right) \mathrm{d}s\\
&=  \sum_{n=0}^\infty f_{n,0} (1-x) \tilde{P}_{n}^{(1,0)}(x)
\end{align*}
Labeling the $(1-x)$ as a weight term and referring to what remains as operator $Q_y$, in reference of it being an integration with respect to $y$, aligns our notation with that in  \cite{gutleb_sparse_2019}. Using reflection operators this can be adapted for the more standard $l(x)=x$ case, see \cite{gutleb_sparse_2019}, which means that considering the kernel $K(x,y)$ means looking at $K(1-x,y)$ on this domain. This operator $Q_y$ acts on a function expanded in the Proriol polynomials with parameters $(0,0,0)$ on $T^2$ and as seen above has the form
\begin{equation*}
\mathrm{Q}_y =   \left(\begin{array}{cccccccccc}
\cline{1-1}
\multicolumn{1}{|c|}{1} &  &  &  &  &  &  &  &  &  \\ \cline{1-3}
\multicolumn{1}{c|}{} & 1 & \multicolumn{1}{c|}{0} &  &  &  &  &  &  &  \\ \cline{2-6}
 &  & \multicolumn{1}{c|}{} & 1 & 0 & \multicolumn{1}{c|}{0} &  &  &  &  \\ \cline{4-10} 
 &  &  &  &  & \multicolumn{1}{c|}{} & \ddots & \ddots & \ddots & \multicolumn{1}{c|}{\ddots} \\ \cline{7-10} 
\end{array}\right).
\end{equation*}
We may account for the as of now omitted weight term $(1-x)$ by using a direct multiplication with Jacobi operators but for reasons of efficiency and due to the need to reflect when $l(x)=x$ is better performed using a bidiagonal lowering operator followed by a diagonal reflection and finally a bidiagonal raising  operator. The discussion so far explains the form of the operator for linear Volterra integral equations of second kind in Algorithm  \ref{alg:linear2} being $\left(\mathbb{1}-\mathrm{S}_{(0,0)}^{(1,0)} \mathrm{R} \mathrm{L}_{(1,0)}^{(0,0)} \mathrm{V}_K \right)$. Equations of first kind are somewhat more subtle and we thus omit discussion of further details, referring instead to the original linear method derivations and proofs in \cite{gutleb_sparse_2019}. We have assumed above that the function may be expanded in the Proriol polynomials but we can make use of additional sparsity structures when instead thinking of $f_{n,k}$ as the extension of univariate function coefficients $f_n$ extended to the triangle domain via the expansion operator
\begin{equation*}
    \mathbf{P}(x,y)^\mathsf{T} \mathbf{f}_\triangle = \mathbf{P}(x,y)^\mathsf{T}  \mathrm{E}_y \mathbf{f}.
\end{equation*}
Choosing the respectively optimal bases $\tilde{\mathbf{P}}^{(1,0)}(x)$ and $\mathbf{P}^{(0,0,0)}(x,y)$ for this purpose results in the extension operator found in \cite{gutleb_sparse_2019}, which when multiplied with the integration from $0$ to $1-x$ operator above results in the following diagonal operator
\begin{equation*}\label{eq:diagonal}
(\mathrm{Q}_y\mathrm{E}_y)_{n,n} =(\mathrm{D}_y)_{n,n}=\frac{(-1)^{n+1}}{n}.
\end{equation*}
Via certain quasi-commutativity properties of the above-discussed operators and the Jacobi operators on the triangle domain and using diagonal reflection operators appropriately one may iteratively build the full Volterra integral operator for a general kernel via the efficient operator-valued Clenshaw algorithm for general kernel linear Volterra integral equations introduced in \cite{gutleb_sparse_2019}. We will refer to this Volterra operator \emph{without} the weight $(1-x)$ as $V_K$ in the following sections and assume it is computed using the methods outlined here and detailed in \cite{gutleb_sparse_2019}. The weight is accounted for by using appropriate basis shifts or multiplication as detailed in the algorithm steps.

\begin{algorithm}
\begin{align*}\int_0^{x} K(x,y) u(y) \mathrm{d}y = g(x).
\end{align*}
\hrulefill
\begin{enumerate}
\item Expand $q(x) = \frac{g(1-x)}{1-x}$ in $\mathbf{\tilde{P}}^{(1,0)}(x)$.
\item Generate $\mathrm{V}_K$ recursively via an operator-valued Clenshaw algorithm for the flipped kernel $K(1-x,y)$.
\item Solve the linear system $\mathrm{V}_K \mathbf{u} = \mathbf{q}$ for $\mathbf{u}$.
\item The approximate solution is $\mathbf{\tilde{P}}^{(1,0)}{}(x)^\mathsf{T} \mathbf{u}$.
\end{enumerate}
\caption{Linear Volterra integral eq. of first kind \cite{gutleb_sparse_2019}}\label{alg:linear1}
\end{algorithm}
\begin{algorithm}
\caption{Linear Volterra integral eq. of second kind \cite{gutleb_sparse_2019}}\label{alg:linear2}
\begin{align*}u(x) = g(x) + \int_0^{x} K(x,y) u(y) \mathrm{d}y.
\end{align*}
\hrulefill
\begin{enumerate}
\item Expand $g(x)$ in $\mathbf{\tilde{P}}^{(1,0)}(x)$.
\item Generate $\mathrm{V}_K$ recursively via an operator-valued Clenshaw algorithm for the flipped kernel $K(1-x,y)$.
\item Solve the linear system $\left(\mathbb{1}-\mathrm{S}_{(0,0)}^{(1,0)} \mathrm{R} \mathrm{L}_{(1,0)}^{(0,0)} \mathrm{V}_K \right) \mathbf{u} = \mathbf{g}$ for $\mathbf{u}.$
\item The approximate solution is $\mathbf{\tilde{P}}^{(1,0)}{}(x)^\mathsf{T} \mathbf{u}$.
\end{enumerate}
\end{algorithm}

\section{Extension of the linear case sparse spectral method to integro-differential equations} \label{sec:integrodifferential}
Volterra integro-differential equations (VIDEs) are named such because the unknown appears in the equation under the action of both a Volterra integral and a derivative operator. In this section we will consider linear VIDEs of the following generic form:
\begin{align}\label{eq:VIDEgeneralcase}
\sum_{m=0}^M \lambda_m \frac{\mathrm{d}^m}{\mathrm{d}x^m}u(x) = g + \mathcal{V}_K u,
\end{align}
with constants $\lambda_m$ and $M \in \mathbb{N}$. Within the context of the present spectral method the integral operator is the Volterra operator from section \ref{sec:volterrabanded}, which as we saw maps a coefficient vector of a function in the $\mathbf{\tilde{P}}^{(1,0)}(x)$ basis to the solution in the same basis. Consistent basis considerations, specficically the Jacobi parameters of our chosen basis, are crucial when developing a solution algorithm for integro-differential equations. As noted in section \ref{sec:functionapproximationuni}, there is an abundance of useful structure in the Jacobi polynomials which among other things allows us to take derivatives by shifting the basis parameters as in (\ref{eq:derivativeop1}-\ref{eq:derivativeop3}). Applying a derivative operator is the same as applying a parameter-scaled raising operator and thus incurs a basis change, which needs to be accounted for in the other operators. We choose the integro-differential equation of second order
\begin{equation*}
\frac{\mathrm{d}^2}{\mathrm{d}x^2} u(x) = g(x)+\int_0^x (x-y)u(y)\mathrm{d}y.
\end{equation*}
as an example to illustrate this. To consistently obtain a solution from an extension to the above linear method it does not suffice to simply replace the second order derivative operator with the appropriate Jacobi polynomial basis derivative operator. Instead, due to the incurred basis shift, an additional conversion or shift operator must be applied to the Volterra operator as well. Starting from the $\mathbf{\tilde{P}}^{(1,0)}(x)$ basis in which we obtain our solution $\mathbf{u}$, the second derivative operator carries us into the basis $\mathbf{\tilde{P}}^{(3,2)}(x)$, meaning that, taking note of the steps in Algorithm \ref{alg:linear2}, the appropriate operator form of the above second order example equation is
\begin{equation*}
\mathbf{\tilde{P}}^{(3,2)}{}^\mathsf{T}\left(\mathcal{D}_2^{(1,0)}-\mathrm{S}_{(1,0)}^{(3,2)}\mathrm{S}_{(0,0)}^{(1,0)} \mathrm{R} \mathrm{L}_{(1,0)}^{(0,0)} \mathrm{V}_K\right)\mathbf{u}_{(1,0)} = \mathbf{\tilde{P}}^{(3,2)}{}^\mathsf{T} \mathbf{g}_{(3,2)}.
\end{equation*}
We may collapse the compatible conversion operators down into a single one to obtain the slightly simpler
\begin{equation}\label{eq:example2ndord}
\mathbf{\tilde{P}}^{(3,2)}{}^\mathsf{T}\left(\mathcal{D}_2^{(1,0)}-\mathrm{S}_{(0,0)}^{(3,2)} \mathrm{R} \mathrm{L}_{(1,0)}^{(0,0)} \mathrm{V}_K\right)\mathbf{u}_{(1,0)} = \mathbf{\tilde{P}}^{(3,2)}{}^\mathsf{T} \mathbf{g}_{(3,2)}.
\end{equation}
Note that $g(x)$ must be expanded in the $ \mathbf{\tilde{P}}^{(3,2)}(x)$ basis instead of the $ \mathbf{\tilde{P}}^{(1,0)}(x)$ basis or converted into said basis using the above-defined basis shift operators. This is for consistency reasons as the operators acting on $\mathbf{u}_{(1,0)}$ shifting the basis from $\mathbf{\tilde{P}}^{(1,0)}(x)$ to $\mathbf{\tilde{P}}^{(3,2)}(x)$ means that the inverse of said operation must act on a function expanded in $\mathbf{\tilde{P}}^{(3,2)}(x)$. This is not an artifact of our choice of the $ \mathbf{\tilde{P}}^{(1,0)}(x)$ basis for our solution: While that basis is particularly well-suited for Volterra integral equations as it results in a far more efficient kernel computation \cite{gutleb_sparse_2019}, the derivative operator in the VIDE will always shift the basis of our solution, so to optimize efficiency $g(x)$ is always initially expanded in the $\mathbf{\tilde{P}}^{(1+M,M)}(x)$ basis where $M$ is the order of the highest appearing derivative operator. Even in the general case with multiple derivative operators of different orders the basis for the solution always remains $\mathbf{\tilde{P}}^{(1,0)}(x)$ (for efficiency of kernel computations) while the highest order derivative operator determines the basis in which $g(x)$ must be expanded along with the shift operators which respectively act on all lower order operators as well as the Volterra integral operator.\\
Attempting to invert the operator on the left-hand side of Equation \eqref{eq:example2ndord} as-is will yield nonsensical results. This should be unsurprising, as the differential equation it corresponds with does not have a unique solution unless initial conditions are supplied as well. Given a Volterra integro-differential equation with highest appearing derivative operator of order $M \in \mathbb{N}$, we will in general require initial conditions for all lower order derivatives to be given, i.e.:
\begin{equation*}
\frac{\mathrm{d}^{m}}{\mathrm{d}x^m}u(0) = c_m, \quad m = 0...M-1,
\end{equation*}
for given constants $c_m$. In the example case of Equation \ref{eq:example2ndord} the values $u(0)$ and $u'(0)$ must be given. In spectral methods such as the one discussed in this paper, boundary or initial conditions are enforced by extending the to-be-inverted operator by appropriate evaluation operators. The relevant Jacobi basis evaluation operators, being functionals, are represented in the coefficient vector and operator language as row vectors, as discussed in section \ref{sec:functionapproximationuni}. For the example in Equation \ref{eq:example2ndord} we thus append the two initial condition evaluations at the top of the operator as follows obtaining the now solvable system
\begin{equation}
\begin{pmatrix}
\mathcal{E}_0\\
\mathcal{E}_0 \mathcal{D}^{(1,0)}\\
\mathcal{D}_2^{(1,0)}-S_{(1,0)}^{(2,3)}V
\end{pmatrix} \mathbf{u}_{(1,0)} =\begin{pmatrix}c_0\\c_1\\
\mathbf{g}_{(3,2)} \end{pmatrix},
\end{equation}
with consistently modified right hand side. Similar procedures have previously been used to solve differential equations, cf. \cite{townsend_automatic_2015,hale_fast_2018}. The discussion in this section in combination with the linear Volterra integral method in Algorithms \ref{alg:linear1} and \ref{alg:linear2} thus provides a recipe for the solution of general linear Volterra integro-differential equations satisfying a sufficient set of initial conditions. We produce the general case method in Algorithm \ref{alg:linearVIDE}. The resulting operator on the left-hand side has filled-in top rows for each initial condition and thus is no longer fully banded but still retains very well-behaved sparsity structure (semi-banded) leading to fast solutions even for high orders of polynomial approximation.
\begin{algorithm}[ht]
\caption{Linear integro-differential Volterra eq. of second kind}\label{alg:linearVIDE}
\begin{align*}\sum_{m=0}^M \lambda_m \frac{\mathrm{d}^m}{\mathrm{d}x^m}u(x) = g(x) + \int_0^{x} K(x,y) u(y) \mathrm{d}y&, \quad \lambda_m \in \mathbb{R}; m,M \in \mathbb{N}\\
\frac{\mathrm{d}^{m}}{\mathrm{d}x^m}u(0) = c_m&, \quad m = 0...M-1, c_m \in \mathbb{R}.
\end{align*}
\hrulefill
\begin{enumerate}
\item Expand $g(x)$ in $\mathbf{\tilde{P}}^{(1+M,M)}(x)$.
\item Generate $\mathrm{V}_K$ recursively via an operator-valued Clenshaw algorithm for the flipped kernel $K(1-x,y)$.
\item Generate the operator $\left(\sum_{m=0}^{M}\lambda_m \mathrm{S}_{(1+m,m)}^{(1+M,M)}\mathcal{D}_m^{(1,0)}-\mathrm{S}_{(0,0)}^{(1+M,M)} \mathrm{R} \mathrm{L}_{(1,0)}^{(0,0)} \mathrm{V}_K\right)$.
\item Append evaluation operators $(\mathcal{E}_0, \mathcal{E}_0 \mathcal{D}^{(1,0)},...)$ to the top row of the operator and corresponding initial conditions $(c_0, c_1, ...)$ to the top of $\mathbf{g}_{(1+M,M)}$.
\item Solve the semi-banded linear system for $\mathbf{u}_{(1,0)}$:\\ $\begin{pmatrix}
\mathcal{E}_0\\
\mathcal{E}_0 \mathcal{D}^{(1,0)}\\
\vdots\\
\mathcal{E}_0 \mathcal{D}^{(M,M-1)}\\
\sum_{m=0}^{M}\lambda_m \mathrm{S}_{(1+m,m)}^{(1+M,M)}\mathcal{D}_m^{(1,0)}-\mathrm{S}_{(0,0)}^{(1+M,M)} \mathrm{R} \mathrm{L}_{(1,0)}^{(0,0)} \mathrm{V}_K
\end{pmatrix} \mathbf{u}_{(1,0)} =\begin{pmatrix}c_0\\c_1\\\vdots\\c_{M-1}\\
\mathbf{g}_{(1+M,M)} \end{pmatrix}$
\item The approximate solution is $\mathbf{\tilde{P}}^{(1,0)}{}(x)^\mathsf{T} \mathbf{u}_{(1,0)}$.
\end{enumerate}
\end{algorithm}

\section{Nonlinear Volterra equations via iterative methods} \label{sec:nonlinear}
In this section we develop an iterative approach for solving nonlinear Volterra integral equations based on the linear case sparse spectral method. Computing solutions to nonlinear Volterra and Fredholm integral equations with iterative methods is not a novel idea in itself, see e.g. \cite{driscoll_automatic_2010}, but typically comes with significant drawbacks, cf. remarks in \cite{atkinson_survey_1992,ezquerro_solving_2011}. The core problem with iterative methods is how rapidly their computational cost scales with the expense of evaluation in each iteration. The slower the rate of convergence and the more expensive the individual evaluation in each step, the less feasible iterative methods become. Conversely, the presented sparse spectral method is very well suited to be used in conjunction with iterative methods as it not only converges exponentially but also keeps evaluation cost comparatively low by making use of operator bandedness in the chosen bases.\\
We will primarily use a simple Newton iteration algorithm without linesearch on the basis of implementations in NLsolve.jl \cite{patrick_kofod_mogensen_nlsolve} for the numerical experiments but in principle many other iterative approaches may be used, resulting in further speed-ups in some cases.\\
The main idea of the extension to the nonlinear case is to notice that given functions $K$, $g$ and $f$ the general nonlinear, second kind Volterra equation
\begin{align*}
u = g + \mathcal{V}_K f(u),
\end{align*}
may be cast into the form of a root-finding problem in function space for the objective function $F(u)$ defined by
\begin{align*}
F(u) := u-\mathcal{V}_Kf(u)-g= 0.
\end{align*}
The initial guess required for iterative approaches is thus made at the level of coefficient vectors, meaning that a guessed column vector representing the solution in the $\mathbf{\tilde{P}}^{(1,0)}(x)$ basis is supplied to the iterative solver. When no convergence automation is used the supplied length of the guess as well as $\mathbf{g}$ determines the maximum polynomial degree and thus the approximation error. The step-by-step method is stated in Algorithm \ref{alg:nonlinearVolterra}.

\begin{algorithm}[ht]
\caption{Nonlinear Volterra integral eq. of second kind}\label{alg:nonlinearVolterra}
\begin{align*}u(x) = g(x) + \int_0^{x} K(x,y) f(y,u(y)) \mathrm{d}y.
\end{align*}
\hrulefill
\begin{enumerate}
\item Expand $g(x)$ in $\mathbf{\tilde{P}}^{(1,0)}(x)$.
\item Generate $\mathrm{V}_K$ recursively via an operator-valued Clenshaw algorithm for the flipped kernel $K(1-x,y)$.
\item Generate the operator $\left(\mathbb{1}-\mathrm{S}_{(0,0)}^{(1,0)} \mathrm{R} \mathrm{L}_{(1,0)}^{(0,0)} \mathrm{V}_K \right)$.
\item Apply a simultaneous root-search (e.g. Newton method) to components of objective function $F(\mathbf{u})=\left(\mathbb{1}-\mathrm{S}_{(0,0)}^{(1,0)} \mathrm{R} \mathrm{L}_{(1,0)}^{(0,0)} \mathrm{V}_K \right)\mathbf{f}(y,\mathbf{u})-\mathbf{g}$.
\item The approximate solution is the obtained root $\mathbf{\tilde{P}}^{(1,0)}{}(x)^\mathsf{T} \mathbf{u}$.
\end{enumerate}
\end{algorithm}

\section{Nonlinear integro-differential Volterra equations}\label{sec:nonlinearintegrodiff}
We can straightforwardly combine considerations in sections \ref{sec:integrodifferential} and \ref{sec:nonlinear} to obtain a sparse spectral method suitable for solving Volterra equations featuring both derivative operators and Volterra integral operators with nonlinearities. A very (but not exhaustively) general case of such an equation of second kind is
\begin{align}\label{eq:nonlinearandintegrodiff}
\sum_{k=0}^m \lambda_k \frac{\mathrm{d}^k}{\mathrm{dx}^k}u(x) = g + \mathcal{V}_K f(u).
\end{align}
For brevity we only address equations of the form in Equation \eqref{eq:nonlinearandintegrodiff} but the methodology outlined in this paper is applicable for a much broader class of problems. The full step-by-step method is stated in Algorithm \ref{alg:nonlinearVIDE}.
\begin{algorithm}
\caption{Nonlinear integro-differential Volterra eq. of second kind}\label{alg:nonlinearVIDE}
\begin{align*}\sum_{m=0}^M \lambda_m \frac{\mathrm{d}^m}{\mathrm{d}x^m}u(x) = g(x) + \int_0^{x} K(x,y) f(y,u(y)) \mathrm{d}y&, \quad \lambda_m \in \mathbb{R}; m,M \in \mathbb{N}\\
\frac{\mathrm{d}^{m}}{\mathrm{d}x^m}u(0) = c_m&, \quad m = 0...M-1.
\end{align*}
\hrulefill
\begin{enumerate}
\item Expand $g(x)$ in $\mathbf{\tilde{P}}^{(1+M,M)}(x)$.
\item Generate $\mathrm{V}_K$ recursively via an operator-valued Clenshaw algorithm for the flipped kernel $K(1-x,y)$.
\item Generate the operator $\left(\sum_{m=0}^{M}\lambda_m \mathrm{S}_{(1+m,m)}^{(1+M,M)}\mathcal{D}_m^{(1,0)}-\mathrm{S}_{(0,0)}^{(1+M,M)} \mathrm{R} \mathrm{L}_{(1,0)}^{(0,0)} \mathrm{V}_K\right)$.
\item Append evaluation operators $(\mathcal{E}_0, \mathcal{E}_0 \mathcal{D}^{(1,0)},...)$ to the top row of the operator and corresponding initial conditions $(c_0, c_1, ...)$ to the top of $\mathbf{g}_{(1+M,M)}$.
\item Apply a simultaneous root-search (e.g. Newton method) to components of objective function $F(\mathbf{u})$ defined by\\ $\begin{pmatrix}
\mathcal{E}_0 \mathbf{u} - c_0\\
\mathcal{E}_0 \mathcal{D}^{(1,0)}\mathbf{u}- c_1\\
\vdots\\
\mathcal{E}_0 \mathcal{D}^{(M,M-1)}\mathbf{u}- c_{M-1}\\
\left(\sum_{m=0}^{M}\lambda_m \mathrm{S}_{(1+m,m)}^{(1+M,M)}\mathcal{D}_m^{(1,0)}-\mathrm{S}_{(0,0)}^{(1+M,M)} \mathrm{R} \mathrm{L}_{(1,0)}^{(0,0)} \mathrm{V}_K\right)\mathbf{f}(y,\mathbf{u}) - \mathbf{g}_{(1+M,M)} \end{pmatrix}
$.
\item The approximate solution is the obtained root $\mathbf{\tilde{P}}^{(1,0)}{}(x)^\mathsf{T} \mathbf{u}$.
\end{enumerate}
\end{algorithm}

\section{Numerical experiments} \label{sec:numexprsection}
Throughout this section we measure errors between analytic solutions $u(x)$ and computed approximate solutions $\mathbf{\tilde{P}}^{(1,0)}(x)^\mathsf{T}\mathbf{u}_{(1,0)}$ in each point of the domain via the infinity norm of the absolute error
\begin{align*}
\| u(x) - \mathbf{\tilde{P}}^{(1,0)}(x)^\mathsf{T}\mathbf{u}_{(1,0)} \|_\infty = \sup_{x\in[0,1]} |u(x)-\mathbf{\tilde{P}}^{(1,0)}(x)^\mathsf{T}\mathbf{u}_{(1,0)}|.
\end{align*}
\subsection{Numerical experiments with linear VIDEs} \label{sec:expintegrodiff}
\subsubsection{Set 1: Second kind, convolution kernels, one derivative operator}
\begin{figure}
    \centering
     \subfloat[]
    {{  \includegraphics[width=5.5cm]{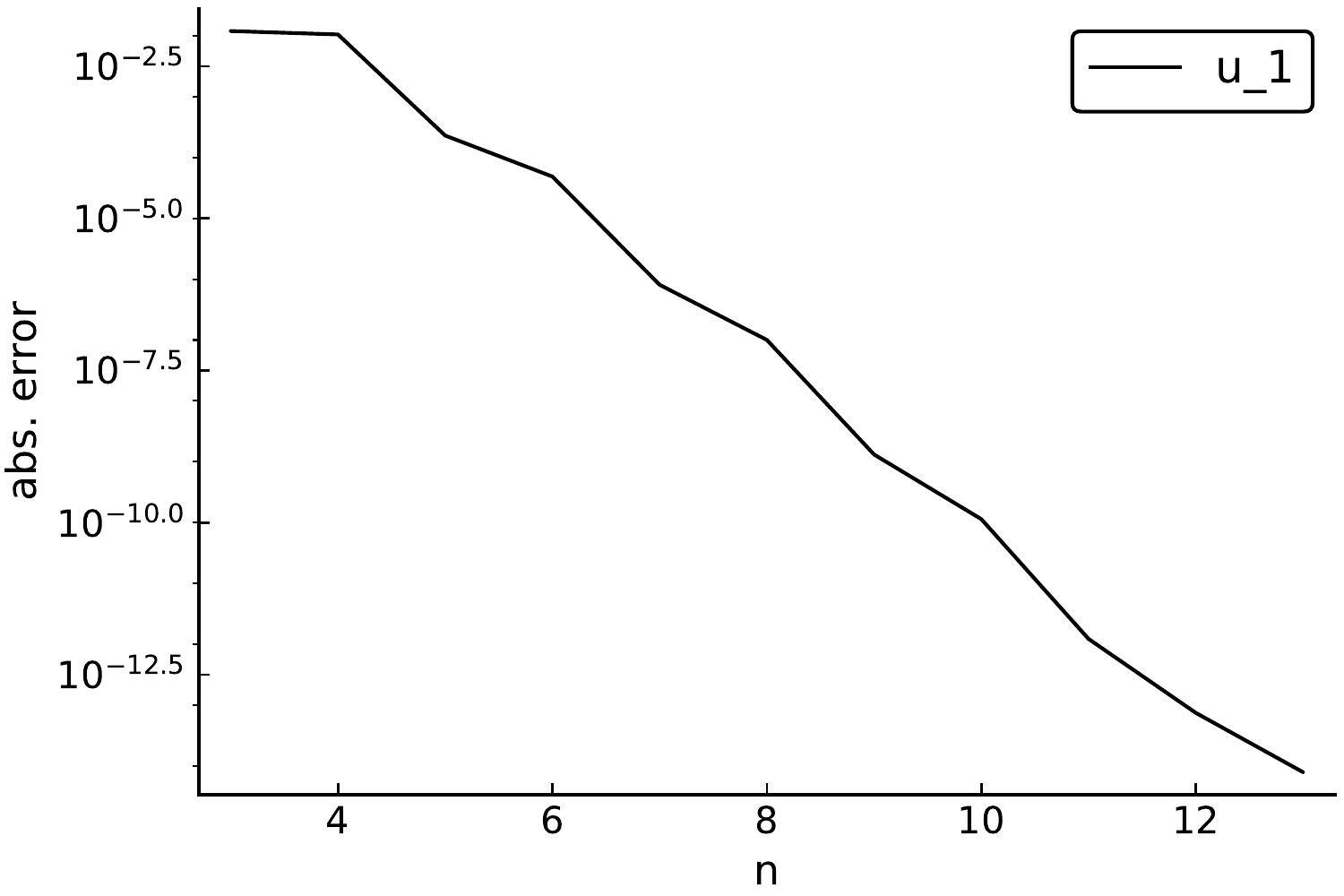} }}
     \subfloat[]
    {{  \includegraphics[width=5.5cm]{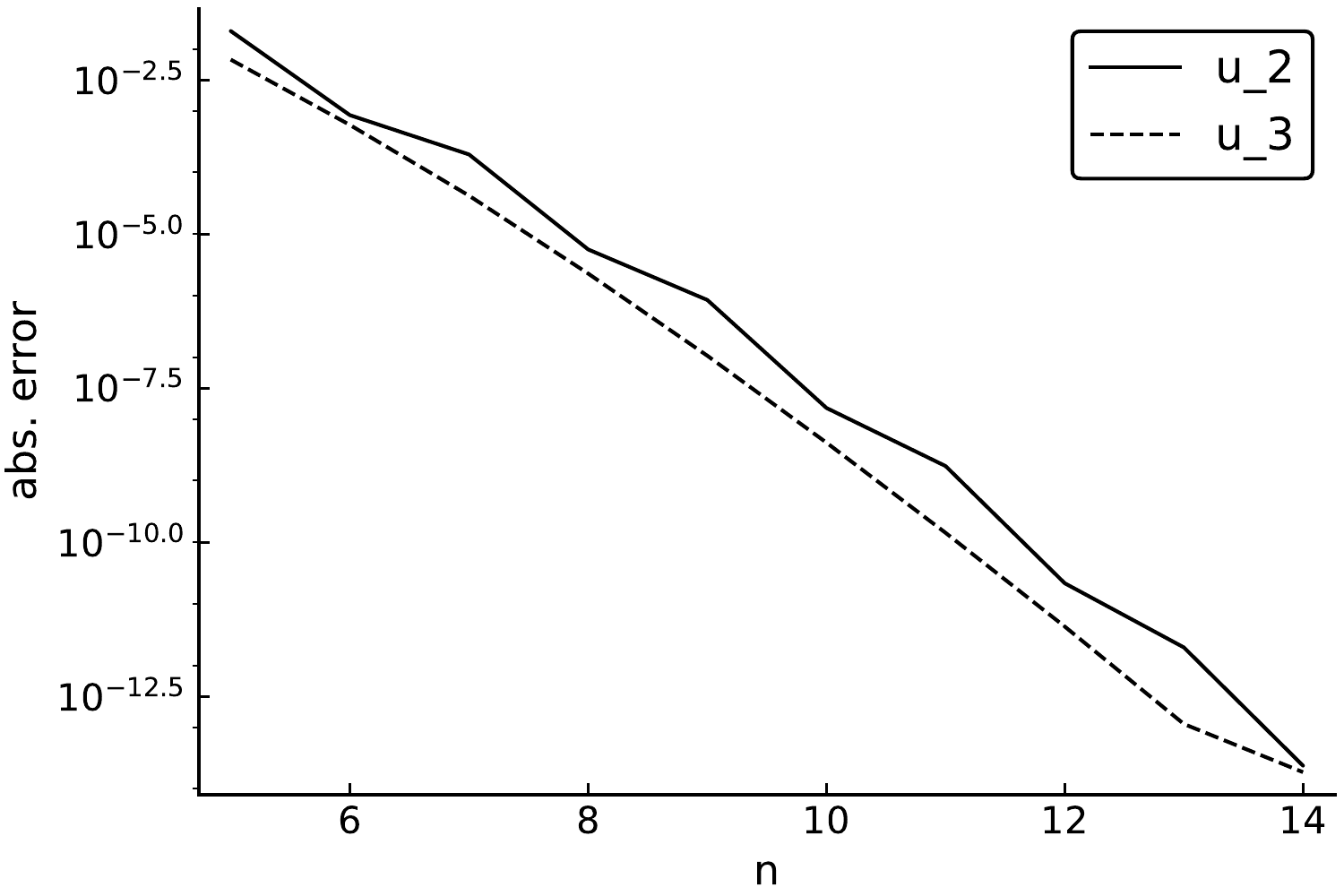} }}
    \caption{Absolute error between analytic and computed solutions for $u_1(x)$, $u_2(x)$ and $u_3(x)$ in equations (\ref{eq:exampleequationss1-1}-\ref{eq:exampleequationss1-3}) for polynomial approximation of order $n$ with initial conditions in (\ref{eq:exampleequationss1initial-1}-\ref{eq:exampleequationss1initial-3}).}%
    \label{fig:set1_errors}%
\end{figure}

As a proof-of-concept we first test the above method on three simple convolution kernel cases with analytically known results: 
\begin{align}\label{eq:exampleequationss1-1}
\frac{\mathrm{d}^2}{\mathrm{d}x^2} u_1(x) = 1+\int_0^x (x-y)u_1(y)\mathrm{d}y,\\
\frac{\mathrm{d}^4}{\mathrm{d}x^4} u_2(x) = -1+x+\int_0^x (y-x) u_2(y) \mathrm{d}y,\\
\frac{\mathrm{d}^3}{\mathrm{d}x^3} u_3(x) = 1+x+\frac{x^2}{2}-\frac{x^4}{4!}+ \int_0^x \frac{(x-y)^2}{2} u_3(y) \mathrm{d}y,\label{eq:exampleequationss1-3}
\end{align}
with initial conditions given by
\begin{align}\label{eq:exampleequationss1initial-1}
u_1(0) = 1,\quad u_1'(0)=0,\\
u_2(0) = -1, \quad u_2'(0) = 1 \quad u_2''(0) = 1 \quad u_2'''(0) = -1,\\
u_3(0)=1, \quad u_3'(0) = 2, \quad u_3''(0)=1.\label{eq:exampleequationss1initial-3}
\end{align}
The following analytic solutions derived respectively via variational iteration, Adomian decomposition and Laplace transform methods are found in \cite{wazwaz_linear_2011}:
\begin{align*}
u_1(x) = \cosh(x),\\
u_2(x) = \sin(x)-\cos(x),\\
u_3(x) = x+e^x.
\end{align*}
We plot the absolute error of the computed solution compared to the analytic solution in semi-logarithmic scale in Figure \ref{fig:set1_errors}, showing exponential convergence to the exact solutions.

\subsection{Set 2: Collocation method for third kind integro-differential equations} \label{sec:comparison_thirdkind}
\begin{figure}
    \centering
     \subfloat[]
    {{  \includegraphics[width=5.5cm]{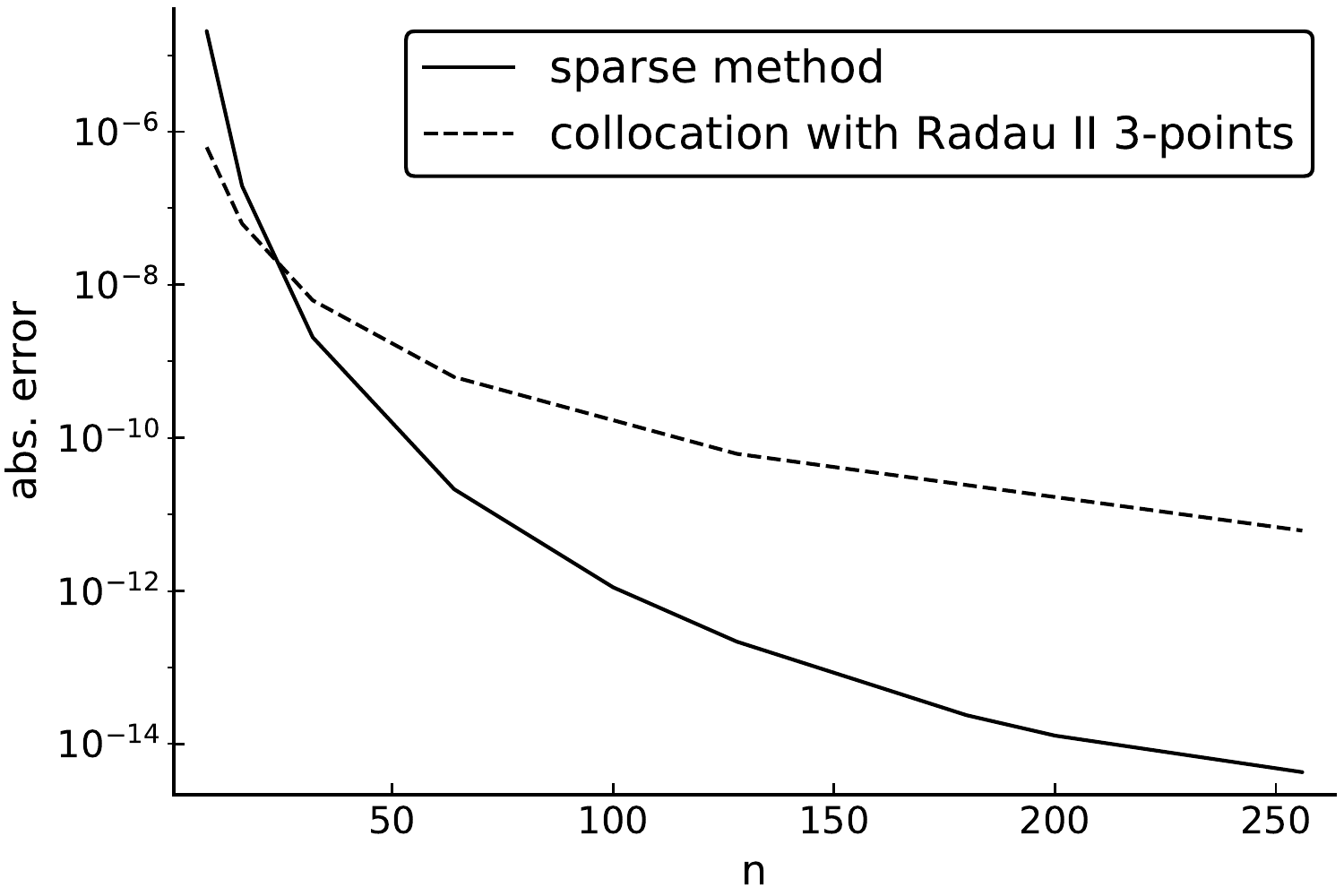} }}
     \subfloat[]
    {{  \includegraphics[width=5.5cm]{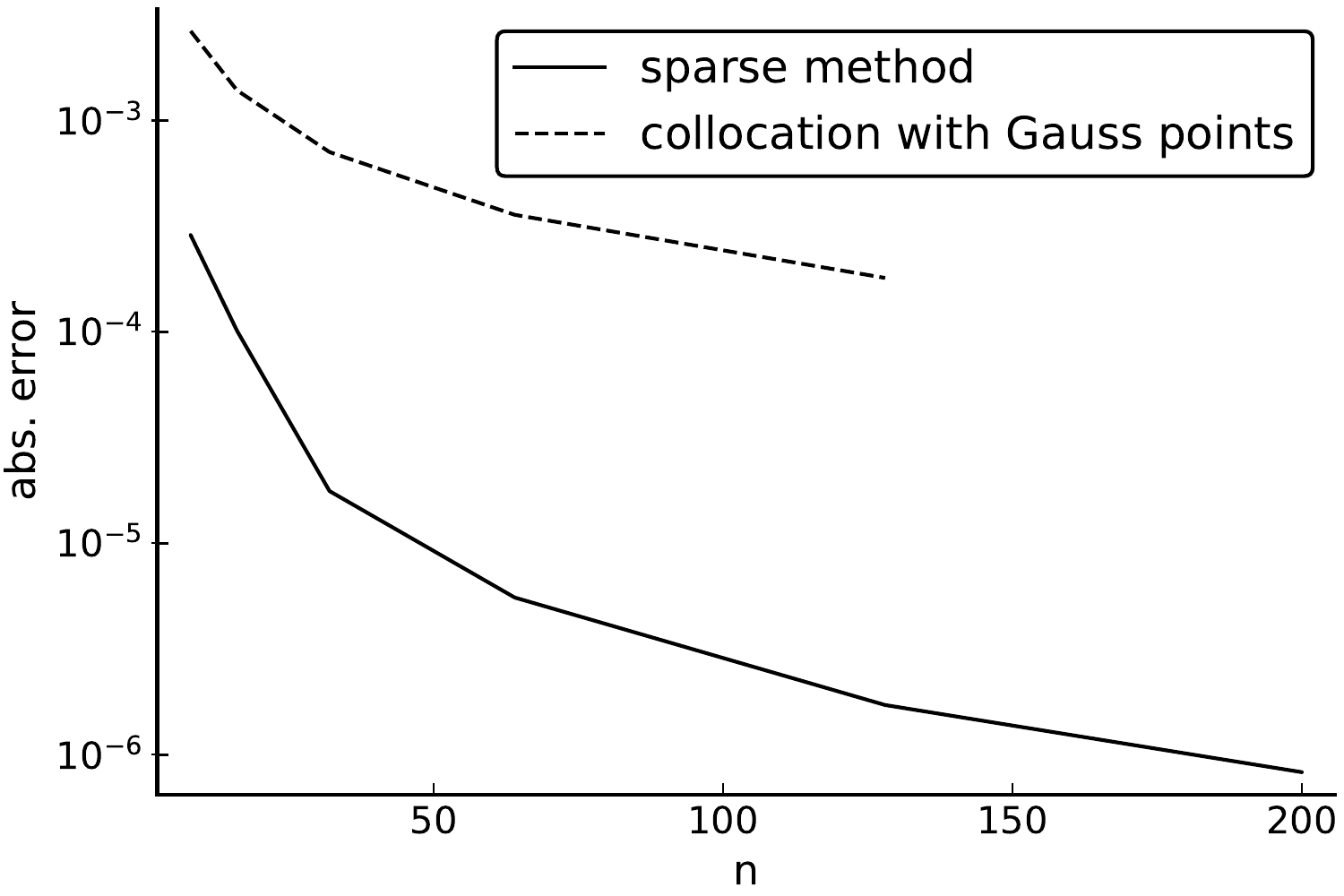} }}
    \caption{(a) shows absolute error between analytic and computed solutions for Equation (\ref{eq:collocation-ex1}) and (b) for Equation (\ref{eq:collocation-ex2}) for polynomial approximations of order $n$. The errors for the collocation method are taken directly from Tables 1 and 2 in \cite{shayanfard_collocation_2019}.}%
    \label{fig:coll_errors}%
\end{figure}

A collocation method is used in \cite{shayanfard_collocation_2019} to solve certain types of third kind integro-differential Volterra equations of form
\begin{equation*}
x^\beta \frac{\mathrm{d}}{\mathrm{d}x}u(x) = x^\beta a(x) u(x) + x^\beta g(x) + \int_0^x K(x,y)u(y)\mathrm{d}y.
\end{equation*}
While we do not explicitly treat third kind equations in this paper, the discussion of first and second kind integro-differential equations in section \ref{sec:integrodifferential} implies an obvious extension to these cases. Shayanfard, Dastjerdi and Ghaini \cite{shayanfard_collocation_2019} discuss two numerical examples and provide a table of error values for differently chosen collocation points. The two numerical experiments with non-convolution kernels are
\begin{align}\label{eq:collocation-ex1}
x^{\frac{2}{3}} u_1'(x) = x^{\frac{2}{3}}\left(\frac{10}{3}x^{\frac{7}{3}}-\frac{3}{16}x^{\frac{14}{3}} \right) + \int_0^x y u_1(y) \mathrm{d}y,\\
x^{\frac{1}{2}}u_2'(x) = \frac{1}{20}xu_2(x)+\frac{9}{2}x^4-\frac{1}{20}x^{\frac{11}{2}}-\frac{1}{6}x^6+\int_0^x y^{\frac{1}{2}}u_2(y)dy,\label{eq:collocation-ex2}
\end{align}
with initial conditions respectively given by
\begin{align}
u_1(0) = 0,\quad
u_2(0) = 0,
\end{align}
and known analytic solutions
\begin{align*}
u_1(x) = x^{\frac{10}{3}},\\
u_2(x) = x^{\frac{9}{2}}.
\end{align*}
In Figure \ref{fig:coll_errors} we compare absolute errors of the results obtained with our sparse spectral approach to the errors obtained with their collocation method (as given in Tables 1 and 2 in \cite{shayanfard_collocation_2019}). The sparse method bandwidth of the operators for these third kind VIDEs is large as they feature additional multiplications with Jacobi operators with poorly approximated rational powers in them, so while our proposed method still performs very well in terms of accuracy, there is not much computation time efficiency to be gained due to sparsity for these two problems.

\subsection{Set 3: Integro-differential equations in Chebfun} \label{sec:comparison_chebfun}
The Chebfun package allows state-of-the-art computations using polynomial approximations and collocation methods in MATLAB \cite{pachon_piecewise-smooth_2010,battles_extension_2004,driscoll_chebop_2008}. An implementation of an automatic collocation method for integral and integro-differential Volterra and Fredholm equations in Chebfun was presented in \cite{driscoll_automatic_2010}. In this section we aim to compare performance of the sparse method compared to the dense collocation method used in Chebfun for problems requiring low and high polynomial orders.
\subsubsection{Low order solutions}
\begin{figure}
    \centering
     \subfloat[]
    {{  \includegraphics[width=5.5cm]{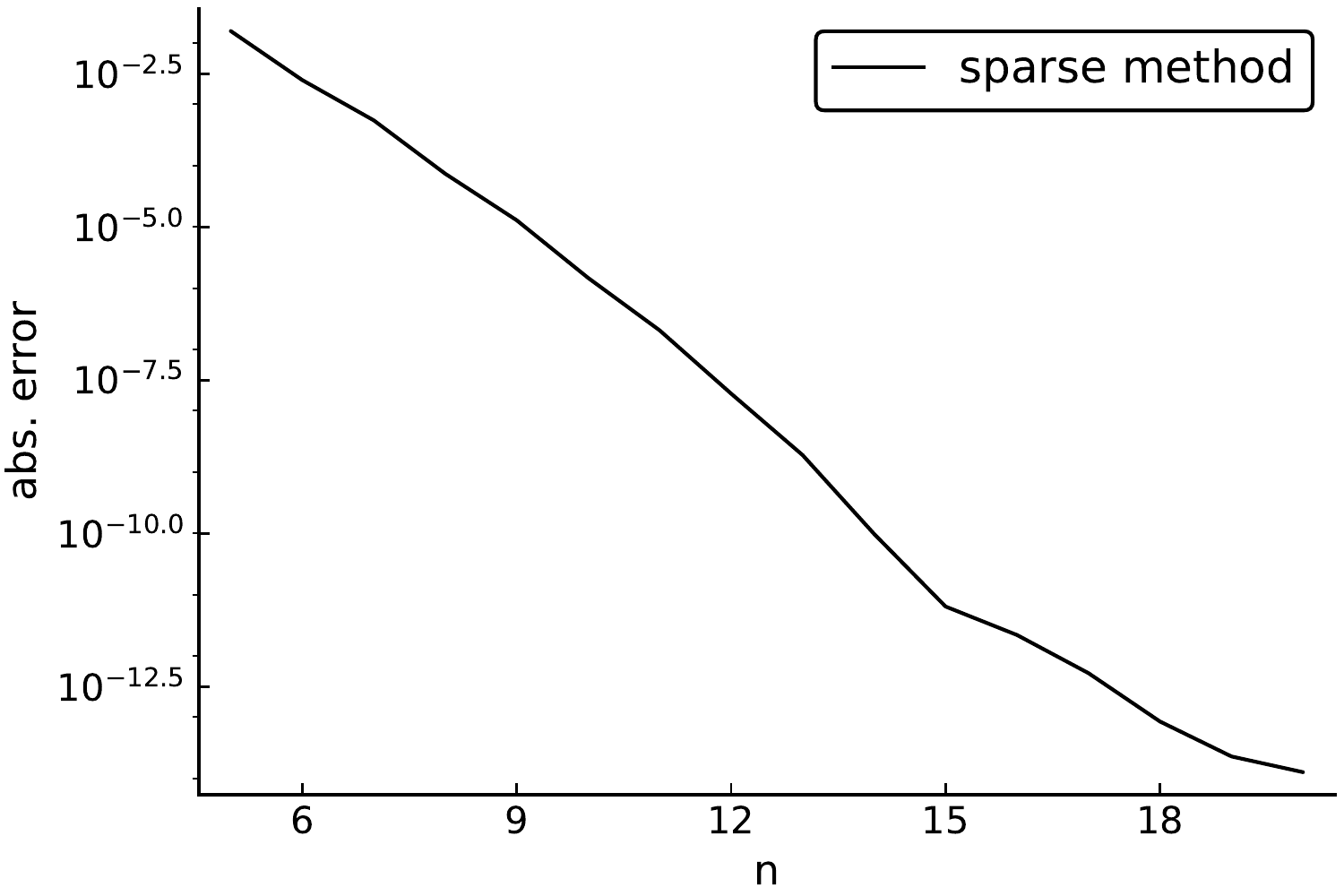} }}
     \subfloat[]
    {{  \includegraphics[width=5.5cm]{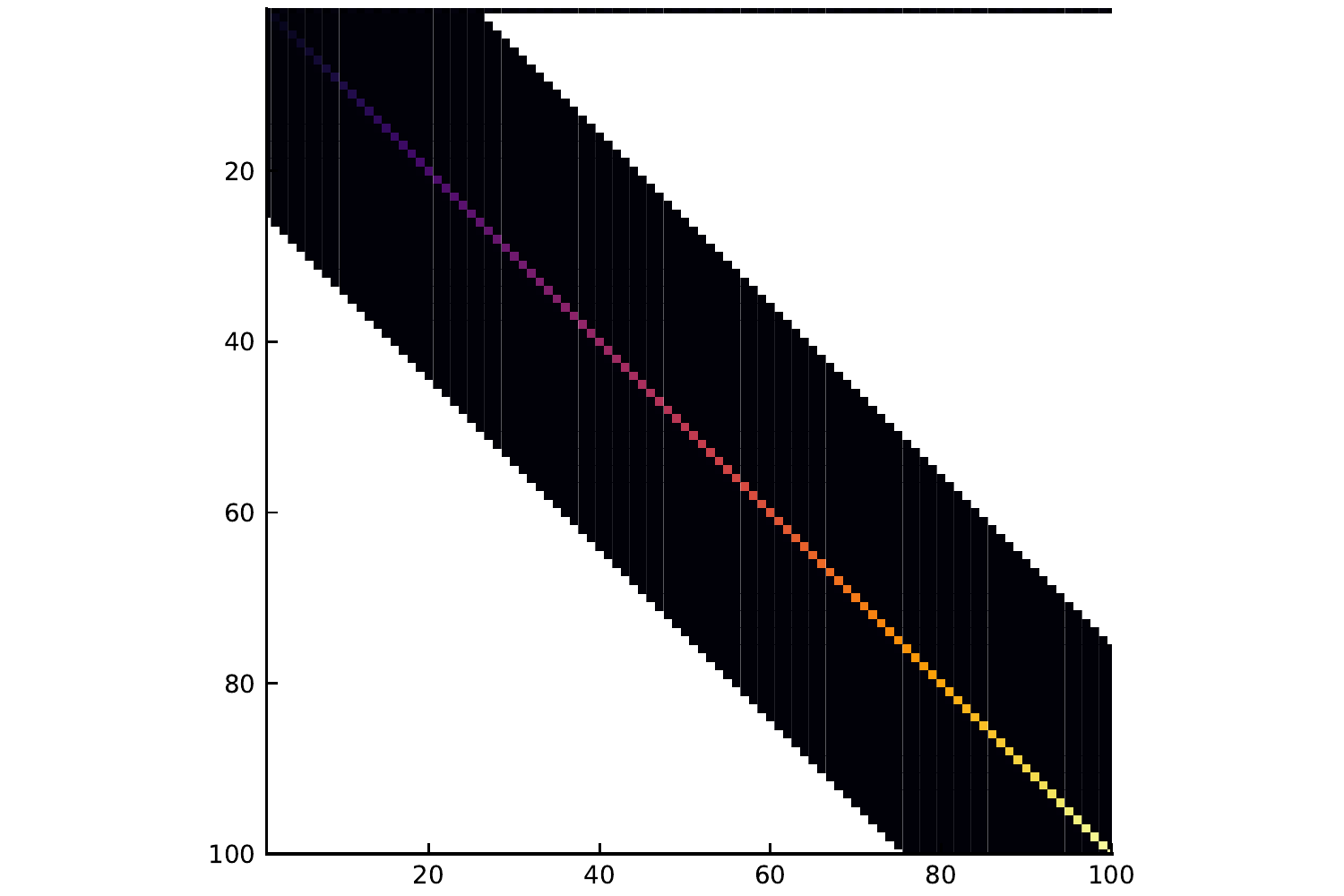} }}
    \caption{(a) shows absolute error between analytic and computed solutions for $u_1$ in equation (\ref{eq:chebfunex-1}) for polynomial approximation of order $n$, (b) shows quasi-bandedness of full sparse method operator for $n=100$.}%
    \label{fig:setChebfun1}%
\end{figure}

The example in this section is given in \cite{driscoll_automatic_2010} and is a non-convolution kernel linear VIDE which previously appeared in a discussion of higher order collocation methods for VIDEs by Brunner \cite{brunner_high-order_1986}. We seek a solution to
\begin{equation}\label{eq:chebfunex-1}
u_1'(x) + u_1(x) = 1+2x+\int_0^x x(1+2x)e^{y(x-y)}u_1(y)\mathrm{d}y,
\end{equation}
with initial condition
\begin{equation}\label{eq:chebfunex-1-initial}
u_1(0) = 1,
\end{equation}
and known analytic solution
\begin{equation*}
u_1(x) = e^{x^2}.
\end{equation*}
In Figure \ref{fig:setChebfun1}(a) we plot the absolute error of the solution obtained via the sparse spectral method with maximal polynomial approximation order $n$. We present a spy plot of the quasi-banded integro-differential operator generated by our sparse method in Figure \ref{fig:setChebfun1}(b). We find that for orders around $n=20$, where machine precision accuracy is within reach, the operator for this problem is still dense and thus the proposed method should realistically only match Chebfun's performance. That a speed-up is nevertheless observed, see Table \ref{tab:chebfun1}, may be explained by language specific differences between MATLAB and Julia, the automatic convergence search which Chebfun performs but was not used for the sparse method or a combination of such factors. Sparsity becomes an important factor for efficiency when treating equations where more complicated solutions are to be expected which require polynomial approximations in the order of hundreds or thousands of coefficients. 
\subsubsection{High order solutions}
\begin{figure}
    \centering
     \subfloat[]
    {{  \includegraphics[width=5.5cm]{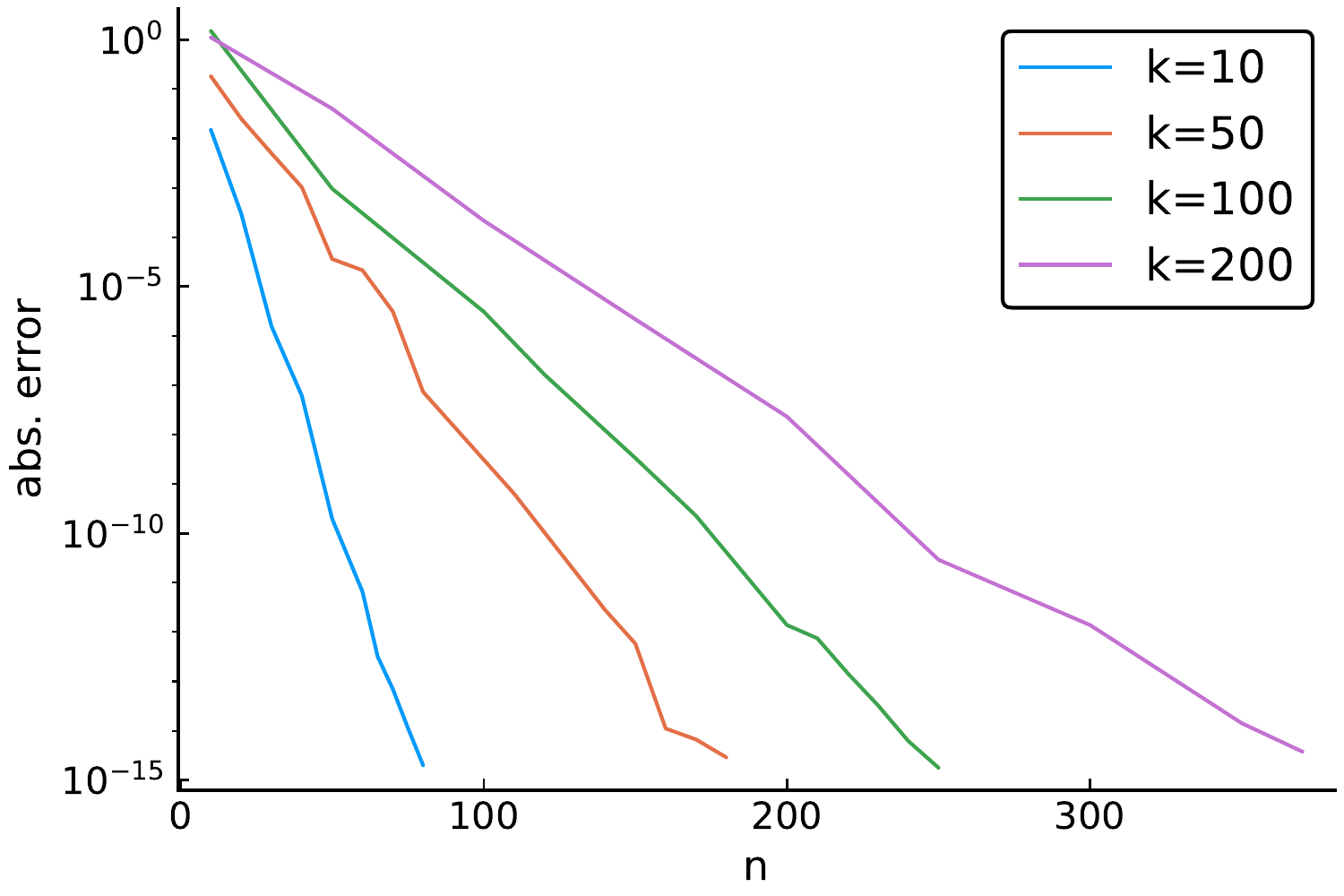} }}
     \subfloat[]
     {{  \includegraphics[width=5.5cm]{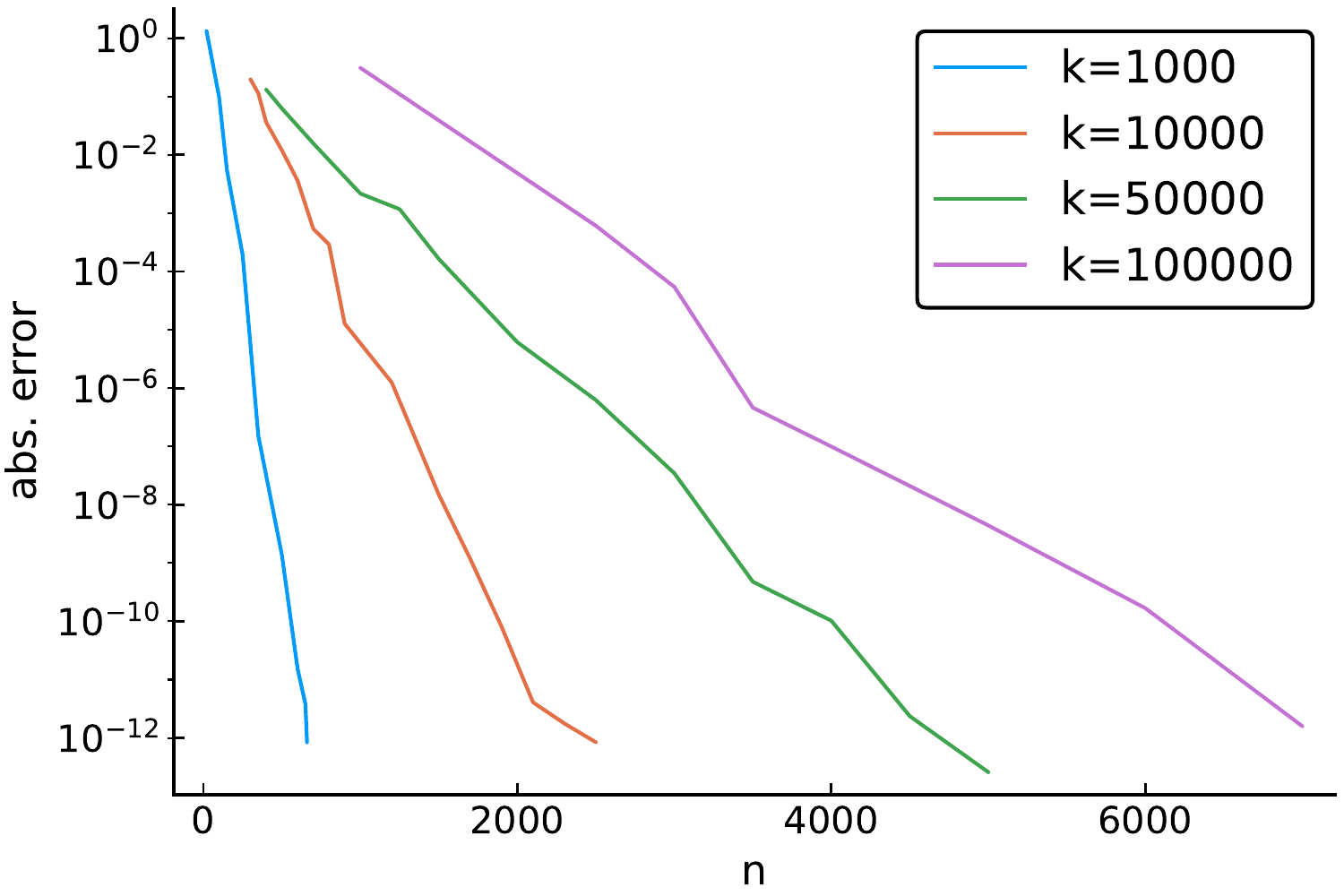} }}
    \caption{(a) and (b) show absolute error between analytic and computed solutions for Equation (\ref{eq:chebfunex-15}) for various values of $k$.}%
    \label{fig:setChebfun15}%
\end{figure}

\begin{figure}
    \centering
         \subfloat[]
    {{  \includegraphics[width=5.5cm]{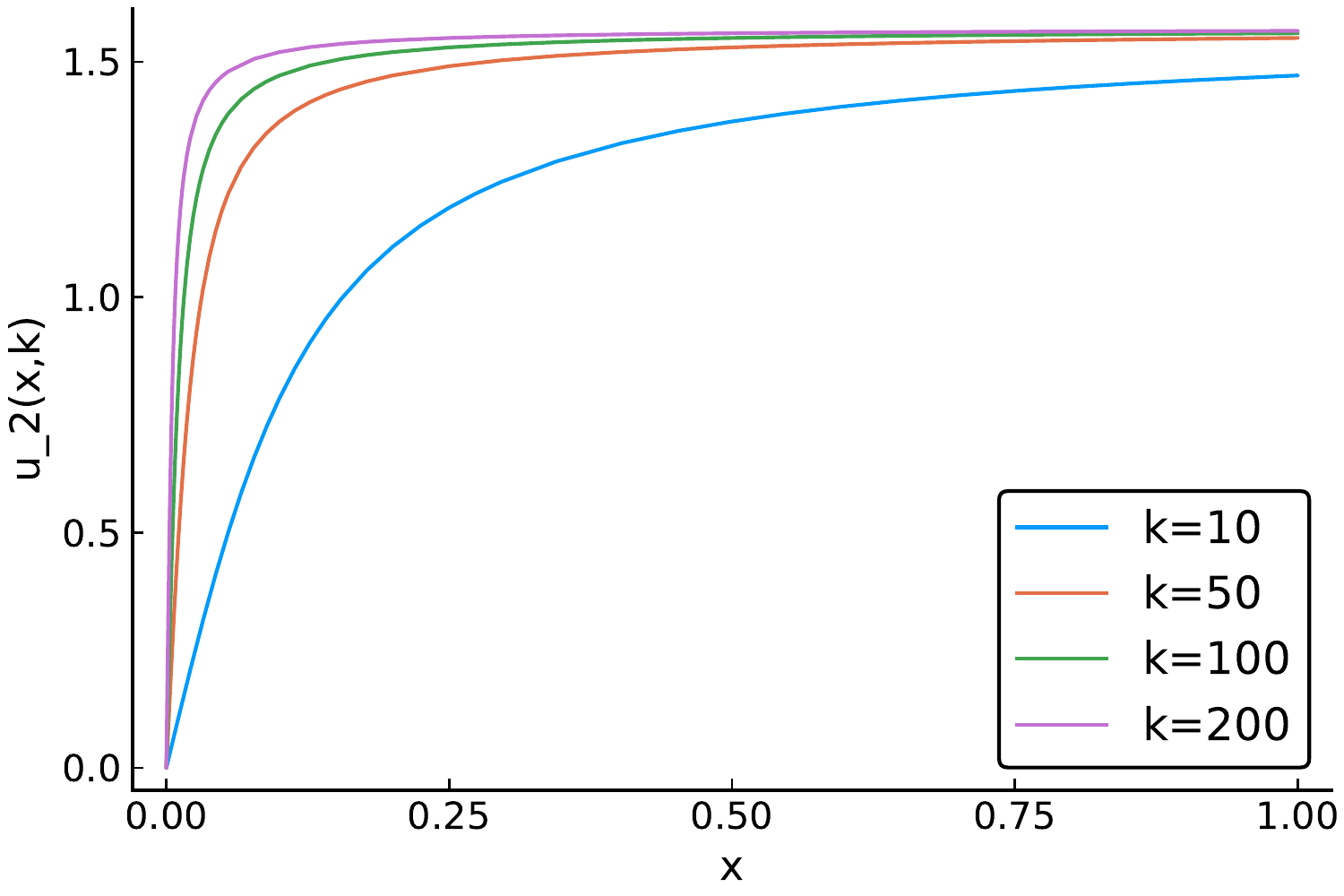} }}
     \subfloat[]
     {{  \includegraphics[width=5.5cm]{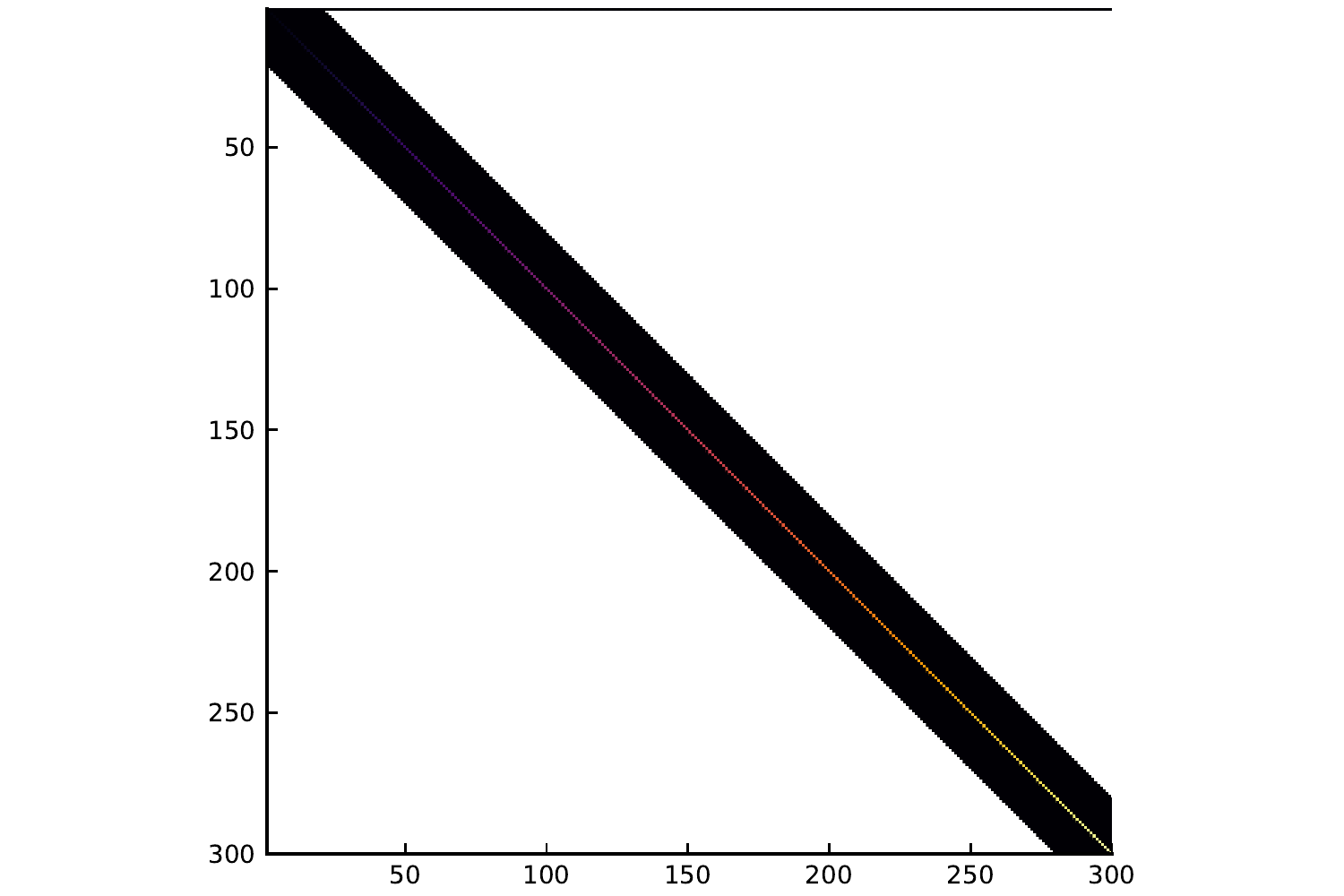} }}
    \caption{(a) shows $u_2(x,k) = \arctan(k x)$ for some small values of $k$, (b) shows quasi-bandedness of full sparse method operator for Equation (\ref{eq:chebfunex-15}) with $n=300$. This banded structure makes computations for very high $n$ not only possible but also fast.}%
    \label{fig:setChebfun15-additionals}%
\end{figure}

For an example which requires a higher $n$ to solve with good accuracy we consider
\begin{equation}\label{eq:chebfunex-15}
u_2'(x,k) = g_2(x,k)+\int_0^x y e^{x^2} u_2(y,k)\mathrm{d}y,
\end{equation}
given initial condition
\begin{equation}\label{eq:chebfunex-15-initial}
u_2(0,k) = 0,
\end{equation}
and right-hand-side function $g_2(x,k)$ defined by
\begin{align*}
g_2(x,k) = \frac{k}{k^2 x^2+1}-\frac{e^{x^2} \arctan(k x)}{2 k^2}+\frac{e^{x^2} x}{2 k}-\frac{1}{2} e^{x^2} x^2 \arctan(k x).
\end{align*}
For all $k \in \mathbb{R}$ the analytic solution to this equation is given by
\begin{equation*}
u_2(x,k) = \arctan(k x).
\end{equation*}
As this approximates a step-like function at $x=0$ for increasing $k$, see Figure \ref{fig:setChebfun15-additionals}(a), it is easy to see why polynomial approximations quickly begin to require high orders. Figure \ref{fig:setChebfun15-additionals}(b) shows a spy plot of the quasi-banded integro-differential operator generated by our sparse method, while Figure \ref{fig:setChebfun15} shows the absolute error of some solutions obtained via the sparse spectral method. The solutions needs to be resolved in relatively high order polynomial approximations, so the bandedness of the operator results in notable performance improvements compared to Chebfun's dense collocation method, see Table \ref{tab:chebfun2}.

\subsection{Set 4: Bessel kernels with highly oscillatory solutions} \label{sec:oscillatorynumexp}
\begin{figure}
    \centering
         \subfloat[]
    {{  \includegraphics[width=5.5cm]{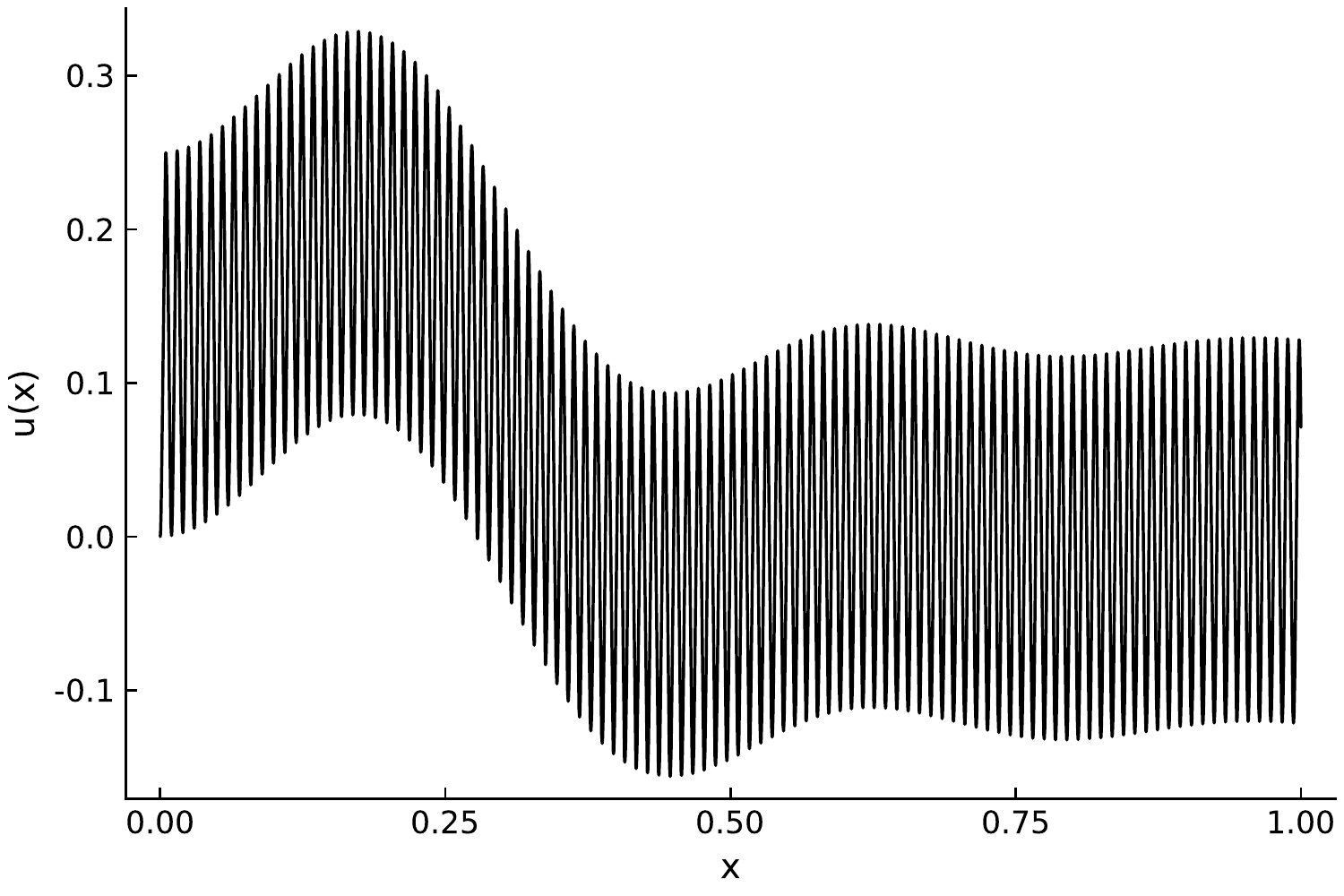} }}
     \subfloat[]
     {{  \includegraphics[width=5.5cm]{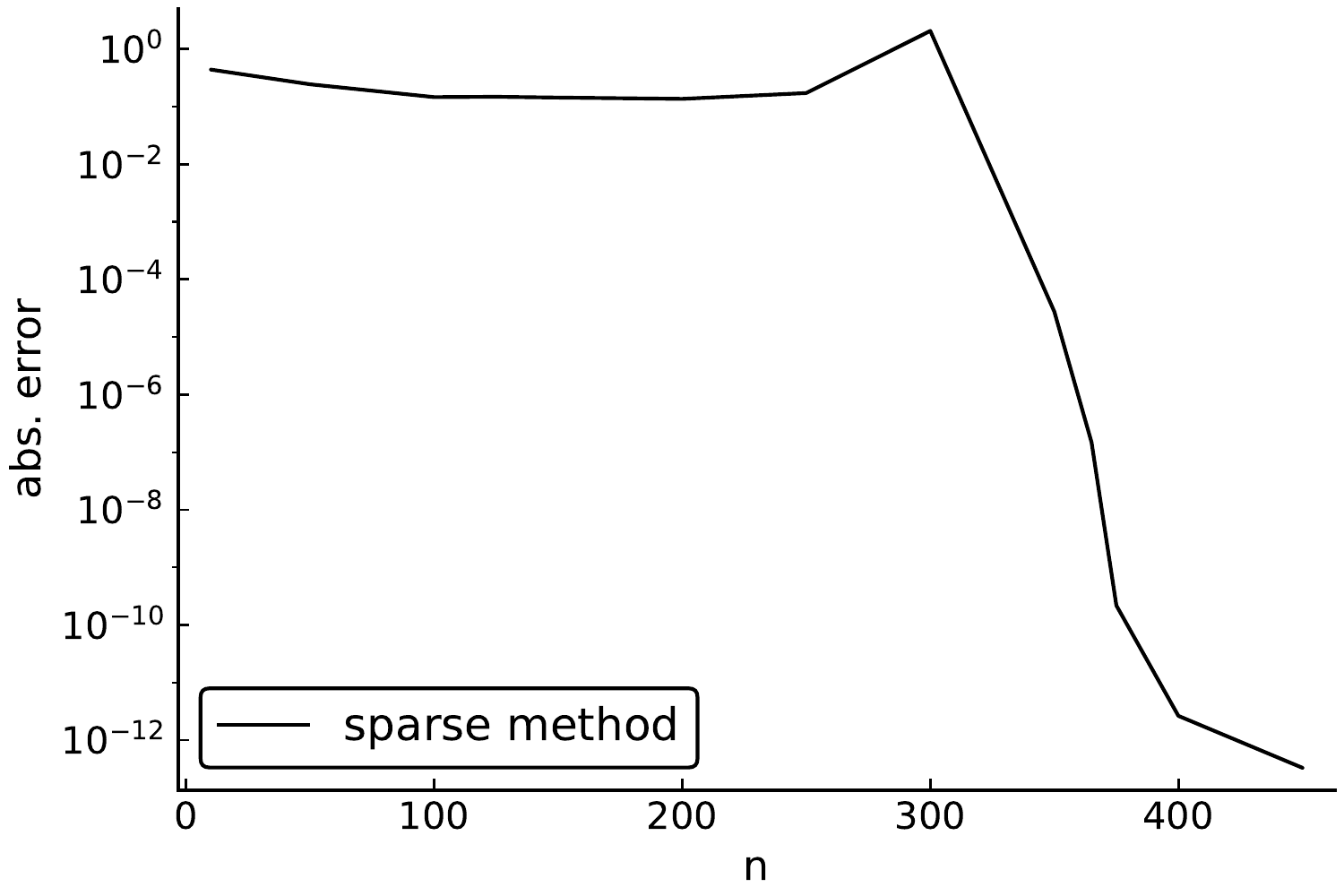} }}
    \caption{(a) shows highly oscillatory solution for Equation \eqref{eq:oscillatoryeq2} for high approximation order $n=2000$, (b) shows error of lower order approximations compared with the $n=2000$ approximation as no analytic solutions are available.}%
    \label{fig:haleoscillatory}%
\end{figure}

Hale \cite{hale_ultraspherical_2019} discusses an integro-differential equation with Bessel function kernel, which appears in scattering applications and potential theory, on the basis of previous work on Bessel kernel Volterra equations by Xiang and Brunner \cite{xiang_efficient_2013}. We use this as the final numerical experiment in this section as it touches on the interesting case of highly oscillatory solutions without known analytic form. Specifically, we discuss the singularly perturbed version of the equation which appears in \cite{hale_ultraspherical_2019} and intentionally makes the problem significantly more oscillatory:
\begin{equation} \label{eq:oscillatoryeq2}
10^{-3} u''(x)+\omega^2 u(x) = g(x,\mu,\nu) - \omega \int_0^x J_\mu (\omega (x-y)) u_2(y) \mathrm{d}y.
\end{equation}
with $g(x,\mu,\nu)$ defined by
\begin{equation*}
g(x,\mu,\nu) = J_{\mu + \nu} (\omega x) + \frac{1}{2x^2}((\nu-1)(\nu-2)J_{\nu-1}(\omega x) + (\nu +1) (\nu+2)J_{\nu+1}(\omega x)),
\end{equation*}
where $J_\mu$ are first kind Bessel functions, $\mu>0$ and $\omega \in \mathbb{R}$. Equation \eqref{eq:oscillatoryeq2} is further supplied with initial conditions
\begin{equation*}
u(0) = u'(0) = 0.
\end{equation*}
To allow comparisons with \cite{hale_ultraspherical_2019} we will consider the example parameters
\begin{align*}
\nu=3, \quad \mu=2, \quad \omega=20.
\end{align*}
Analytic solutions to this equation are not known and convergence comparisons are thus made to high order approximate solutions ($n=2000$) instead. We plot the highly oscillatory solution to this in Figure \ref{fig:haleoscillatory}(a) and the convergence to the $n=2000$ solution in Figure \ref{fig:haleoscillatory}(b). Similarly to results in \cite{hale_ultraspherical_2019} we observe rapid exponential convergence once the polynomial order becomes sufficient to resolve the frequency of the oscillations. Better convergence up to machine precision is possible when using a more sophisticated balancing of approximation orders for the kernel, $g$ and the solution respectively, as opposed to linearly increasing the approximation order of each of them at the same time. This could also be done using an automated convergence algorithm if needed but this example is primarily presented to show the broad range of applicability even for oscillatory problems -- as this particular Bessel kernel is ultimately a convolution kernel, methods which take the additional convolution kernel structure into account, e.g. Hale's method in \cite{hale_ultraspherical_2019}, will generally outperform the general kernel method presented in this paper in accuracy or performance (in particular if operating in low polynomial approximation orders or if they themselves make use of sparsity structure).

\subsection{Numerical experiments with nonlinear equations} \label{sec:expnonlinear}
\subsubsection{Set 1: Power nonlinearity Volterra integral equations}
\begin{figure}
    \centering
     \subfloat[]
    {{  \includegraphics[width=5.5cm]{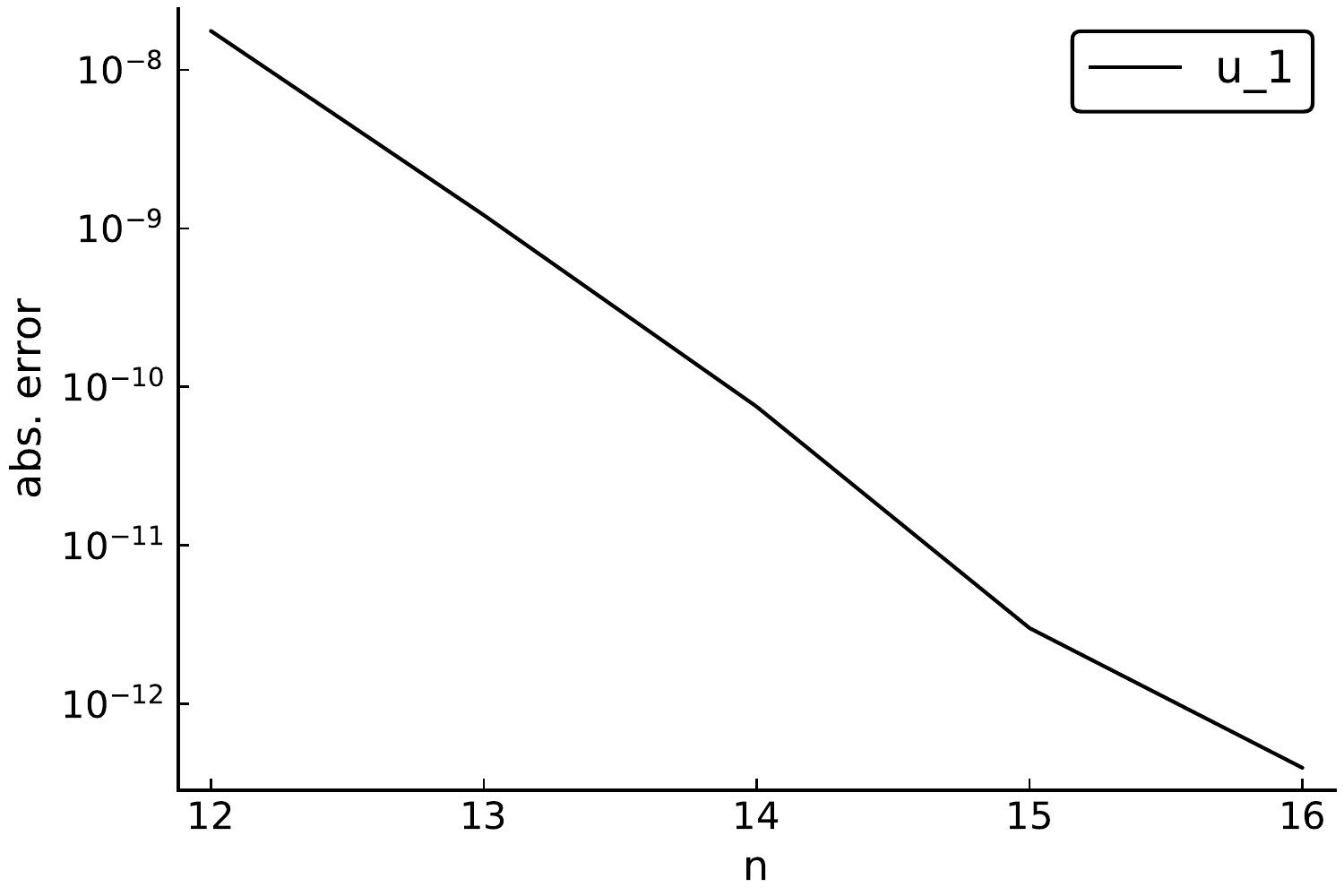} }}
     \subfloat[]
    {{  \includegraphics[width=5.5cm]{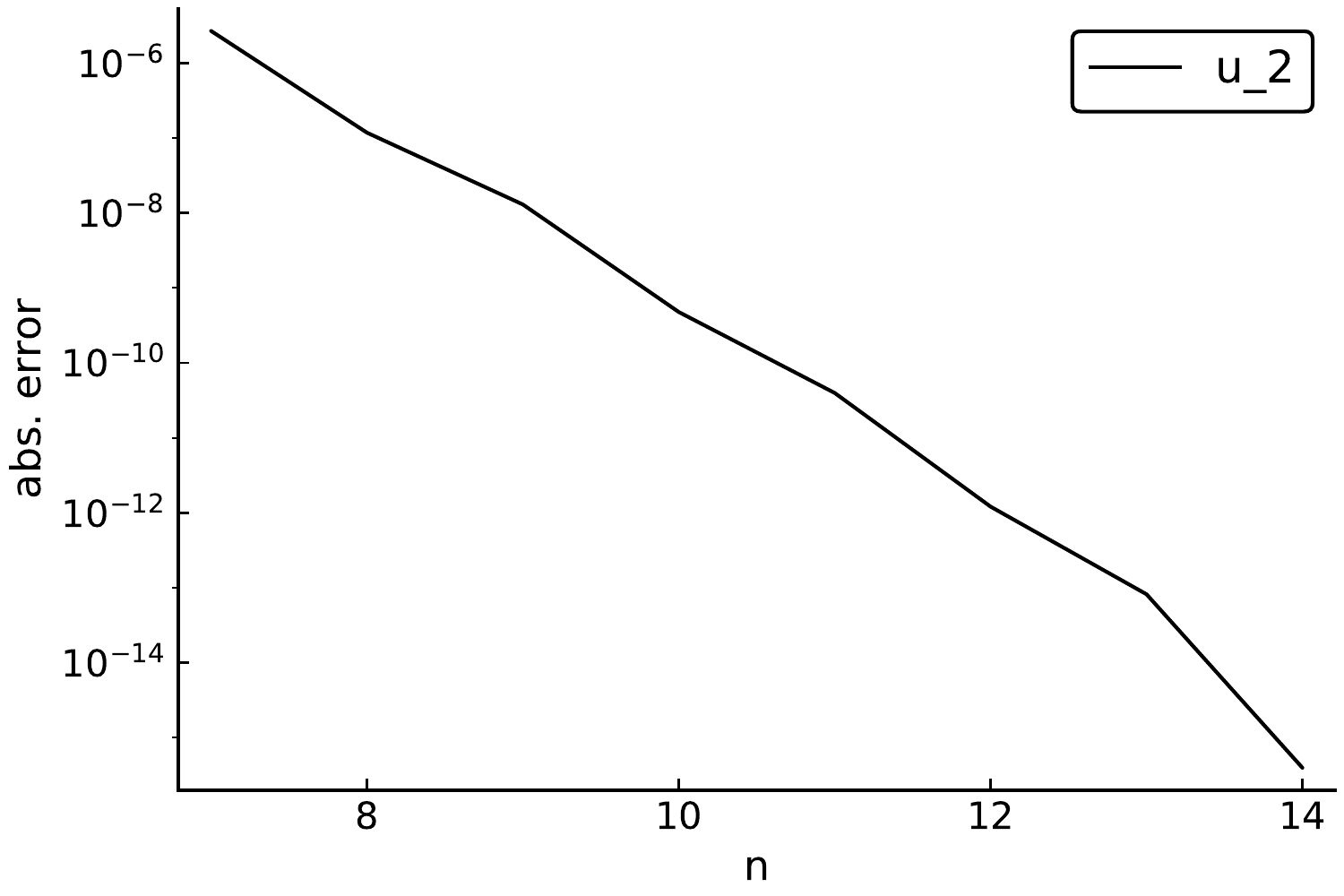} }}
    \caption{Absolute error between analytic and computed solutions for $u_1(x)$ and $u_2(x)$ in equations (\ref{eq:exampleequationssNL1-1}-\ref{eq:exampleequationssNL1-2}) for polynomial approximation of order $n$.}%
    \label{fig:set1NL_errors}%
\end{figure}

The simplest case of nonlinear Volterra integral equations and thus also where most analytic solutions are available for direct comparison is the case of power nonlinearities of the form $f(u)=u^m$ for some positive integer $m$. We thus consider the examples
\begin{align}\label{eq:exampleequationssNL1-1}
u_1(x) = e^x+\frac{x(1-e^{3x})}{3}+\int_0^x x u_1^3(y)\mathrm{d}y,\\
u_2(x) = \sin(x)+\frac{\sin^2(x)}{4}-\frac{x^2}{4}+\int_0^x (x-y) u_2^2(y) \mathrm{d}y,\label{eq:exampleequationssNL1-2}
\end{align}
whose analytic solutions are derived respectively via a Picard type iteration and Adomian decomposition method in \cite{wazwaz_linear_2011}:
\begin{align*}
u_1(x) = e^x,\\
u_2(x) = \sin(x).
\end{align*}
As discussed above and as is true for any iterative method there are now multiple parameters which may be fine-tuned to the problems at hand in order to achieve faster and more precise convergence. To that end we may for example fine-tune the initial guess or the convergence cut-offs. As what can go wrong in a standard application case is of greater interest than what may happen in ideal circumstances, we omit such fine-tuning and instead simply supply a vector of all zeros of length $n$ for Equation \eqref{eq:exampleequationssNL1-1} and a vector of all ones of length $n$ for Equation \eqref{eq:exampleequationssNL1-2}. We plot the maximal absolute errors between true and computed solutions in Figure \ref{fig:set1NL_errors}. We observe exponential convergence as $n$ increases using simple Newton iteration without linesearch.

\subsection{Set 2: Numerical experiments with nonlinear VIDEs} \label{sec:expnonlinearintegrodiff}
\begin{figure}
    \centering
     \subfloat[]
    {{  \includegraphics[width=5.5cm]{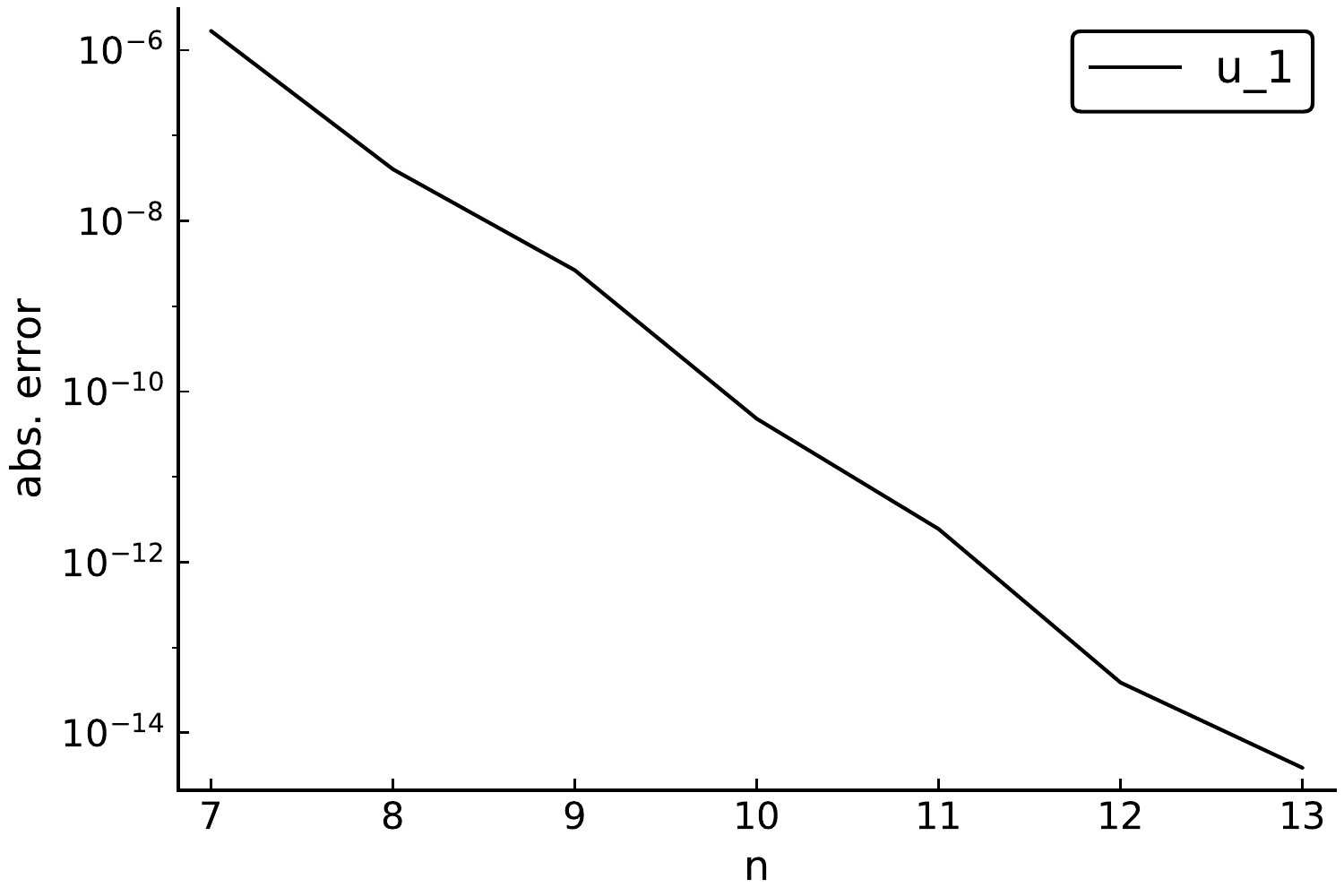} }}
     \subfloat[]
    {{  \includegraphics[width=5.5cm]{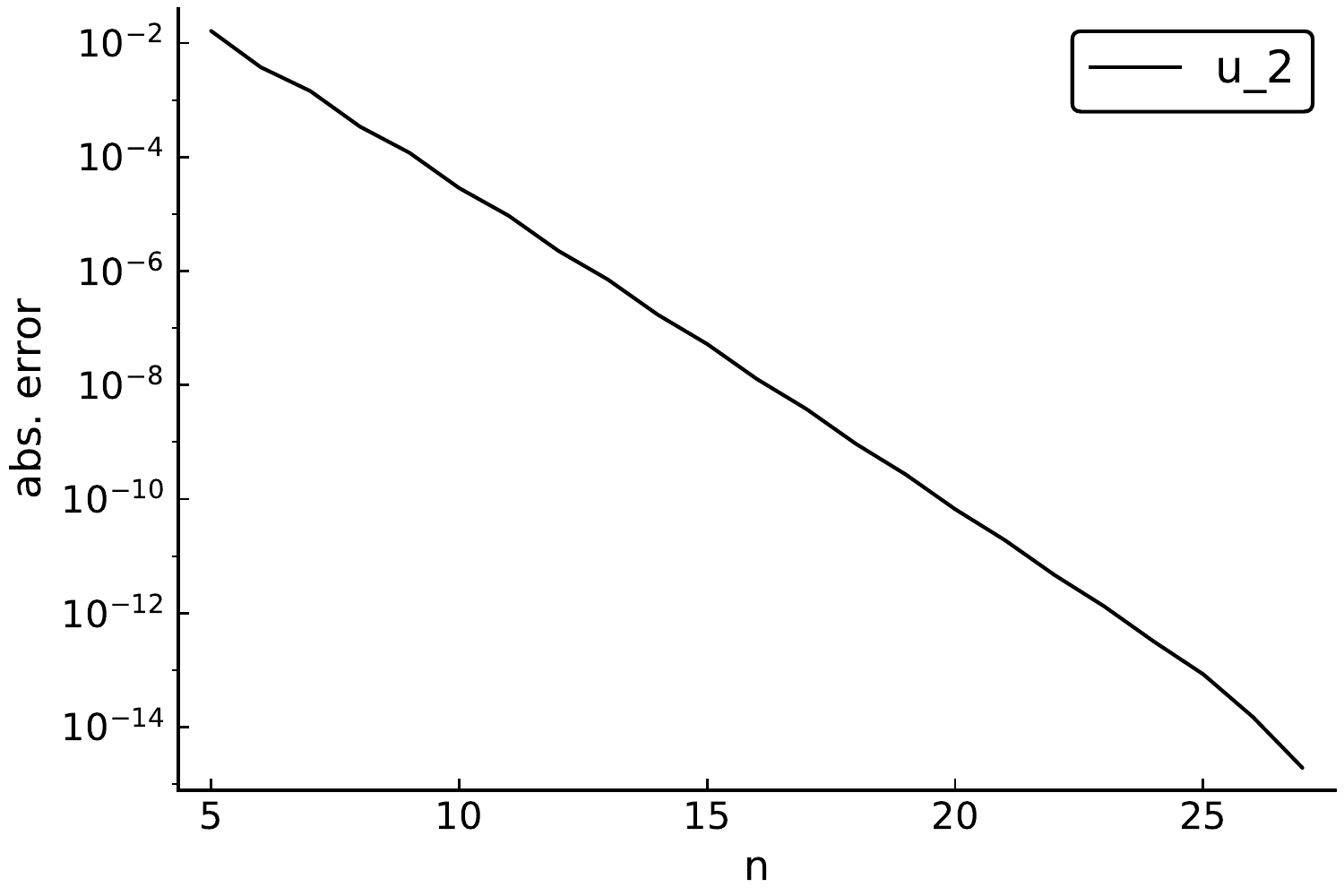} }}
    \caption{Absolute error between analytic and computed solutions for $u_1(x)$ and $u_2(x)$ in equations (\ref{eq:exampleequationsNLInts1-1}-\ref{eq:exampleequationsNLInts1-2}) with initial conditions in (\ref{eq:nonlinearandintegrodiffexpset1inital1}-\ref{eq:nonlinearandintegrodiffexpset1inital2}) for polynomial approximation of order $n$.}%
    \label{fig:set1INTNL_errors}%
\end{figure}

For nonlinear VIDEs we consider:
\begin{align}\label{eq:exampleequationsNLInts1-1}
\frac{\mathrm{d}^2}{\mathrm{d}x^2} u_1(x) = -\frac{5}{3}\sin(x)+\frac{1}{3}\sin(2x)+\int_0^x \cos(x-y)u_1^2(y)\mathrm{d}y,\\
\frac{\mathrm{d}}{\mathrm{d}x} u_2(x) = x+\cos(x)-\tan(x)+\tan^2(x)+ \int_0^x u_2^2(y) \mathrm{d}y,\label{eq:exampleequationsNLInts1-2}
\end{align}
with initial conditions
\begin{align}\label{eq:nonlinearandintegrodiffexpset1inital1}
u_1(0) = 0, \quad u'_1(0) = 1,\\
u_2(0) = 0.\label{eq:nonlinearandintegrodiffexpset1inital2}
\end{align}
Analytic solutions to these equations were derived in \cite{wazwaz_linear_2011} using the variational iteration method:
\begin{align*}
u_1(x) = \sin(x),\\
u_2(x) = \tan(x).
\end{align*}
As in section \ref{sec:expnonlinear} we avoid making educated guesses for the inital guess supplied to the algorithm and merely increase the maximal allowed length of the solution coefficient vector, i.e., the maximal polynomial degree of the computed approximation. The initial guess for Equation (\ref{eq:exampleequationsNLInts1-1}) is a vector of all ones and the initial guess for Equation (\ref{eq:exampleequationsNLInts1-1}) is a vector of all zeros of length $n$ respectively. We plot the maximal absolute errors between analytic and computed solutions in Figure \ref{fig:set1INTNL_errors}. We again observe exponential convergence as $n$ increases using Newton iteration without linesearch.
\begin{table}
\begin{center}
\begin{tabular}{c c  c c c}
\textbf{Method}         &      \textbf{CPU time}       & \textbf{approx. order} &  \textbf{abs. error}    \\ \hline \hline
 Sparse method           & 0.1s                   & 20 & 1.4e-15    \\   
Chebfun            & 0.2s       & 18 (autom.) & 2.7e-15  \\  
\end{tabular}
\end{center}
\caption{Quantitative performance comparison of sparse method and Chebfun for Equation (\ref{eq:chebfunex-1}). The automatically chosen convergence order was used for Chebfun's results. CPU time measured on Intel Core i7-6700T CPU @ $2.80$GHz.}
\label{tab:chebfun1}
\end{table}
\begin{table}
\begin{center}
\begin{tabular}{c c  c c c}
\textbf{k (Sparse)\hspace{1mm} }         &      \textbf{CPU time}       & \textbf{approx. order} &  \textbf{abs. error}  \\ \hline \hline
100            & 0.3s                   & 300 & 6.0e-15     \\   
1000            & 0.4s         & 800  & 2.9e-14   \\  
10000            & 0.8s         & 2000 & 9.3e-14  \\   
50000            & 3.5s   & 5000       & 2.8e-13  \\   
100000 & 8.8s & 8000 & 8.1e-12
\end{tabular}\\
\begin{tabular}{c c  c c c}
\textbf{k (Chebfun)}         &      \textbf{CPU time}       & \textbf{autom. order} &  \textbf{abs. error}     \\ \hline \hline
100            & 1.4s                   & 196 & 4.8e-14 \\   
1000            & 2.1s       & 540 & 6.6e-12   \\  
10000            &  10.7s     & 1477 & 1.2e-09  \\   
50000            & 10.4s        &  2863 & 4.1e-08 
\end{tabular}\\
\end{center}
\caption{Quantitative performance comparison of sparse method and Chebfun for Equation (\ref{eq:chebfunex-15}). The automatically chosen order was used for Chebfun's results, while the sparse method can generate higher accuracy results in less time. For $k=100000$ Chebfun issues an error after approximately $11s$ that it may not have converged. CPU time measured on Intel Core i7-6700T CPU @ $2.80$GHz.}
\label{tab:chebfun2}
\end{table}

\section{Notes on algorithm convergence} \label{sec:analysis}
Convergence of the above-discussed method in the case of nonlinear equations arises as a function of the convergence properties of the root search algorithm that is utilized, combined with the proofs for the respective linear variants. Proofs of convergence for second kind linear Volterra integro-differential equations may be given in a similar fashion to linear Volterra integral equations in \cite{gutleb_sparse_2019} and differential equations in \cite{hale_fast_2018}, the basic observation being that the full to-be-inverted operator is diagonally dominant for well-behaved functions and may be written as a compact perturbation of the identity, thus reducing the problem to standard finite section approximation convergence results, cf. \cite{bottcher_analysis_2006,olver_fast_2013,slevinsky_singular_2017,lintner_generalized_2015}. As seen in the linear case proof for first kind VIEs in \cite{gutleb_sparse_2019} proofs for first kind equations would require a deeper functional analysis approach due to being somewhat ill-conditioned.

\section{Discussion} \label{sec:discussion}
We have presented a competitively fast general kernel sparse spectral method for nonlinear Volterra integro-differential and integral equations which extends linear results in \cite{gutleb_sparse_2019}. The method is notably not reliant on the structure of convolution kernels and applies for general kernels. Furthermore, as it does not rely on low rank approximations it is applicable in more general cases where these approximations fail. It thus combines very broad applicability with high performance and accuracy.\\
One noteworthy drawback of this method is that, although as discussed in the numerical experiments section in \cite{gutleb_sparse_2019} the method may yield sensible results for some types of singular kernels, there are as of now no known guarantees for such cases. That said, the presented method was shown to be convergent and well-behaved with problems that may be well approximated in the specified polynomial bases, which allow for a very general range of kernels.\\
The numerical experiments in this paper serve an illustrative purpose -- in a practical application setting one would choose more sophisticated and efficient root search algorithms than a simple Newton iteration without linesearch and make an educated initial guess for the root search based on background knowledge about the structure of the problem instead of supplying simple zero or one filled coefficient vectors. These points were specifically ignored in this paper to illustrate that such more sophisticated methods are not required to achieve competitive performance and accuracy.

\begin{acknowledgements}
The author would like to thank Sheehan Olver for reading a draft and providing helpful comments.
\end{acknowledgements}

%

\bibliographystyle{spmpsci}      
\bibliography{references.bib}   

\end{document}